\documentclass[12pt]{amsart}
\usepackage{euler,amsmath,amssymb,amscd}
\usepackage{epsfig}

\newtheorem{thm}{Theorem}
\newtheorem{prop}{Proposition}
\newtheorem{coro}{Corollary}
\newtheorem{lemma}{Lemma}

\newtheorem*{thma}{Theorem}

\theoremstyle{definition}

\newtheorem{remark}{Remark}

\newtheorem{ex}{Exercise}
\newtheorem{defn}{Definition}

\newtheorem*{ack}{Acknowledgement}

\def\mbb{\mathbb}
\def\bb{\mathbb}
\def\mcl{\mathcal}

\def\star{\ast}
\def\ten{\otimes}
\def\ex{\times}

\def\tu{\textup}

\def\ndt{\noindent}
\def\lan{\langle}
\def\ran{\rangle}

\def\a{\alpha}
\def\b{\beta}
\def\d{\delta}
\def\D{\Delta}
\def\e{\epsilon}
\def\g{\gamma}

\def\G{\Gamma}
\def\D{\Delta}

\def\la{\lambda}

\def\sm{\sigma}

\def\o{\omega}

\def\P{\mbb P}
\def\Z{\mbb Z}

\def\Q{\mbb Q}

\def\bG{\mathbb G}

\def\inj{\hookrightarrow}

\def\SS{\mcl{O}}
\def\O{\mcl O}

\def\inv{^{-1}}

\def\Proj{\tu{\textbf{Proj}\,}}

\def\mod{/ \! \! /}
\def\hilb{Hilb}
\def\SL{SL}
\def\om2{\omega^{\ten 2}}
\def\Gr{Gr}
\def\Mg{\overline{M}_g}

\def\H{Hilb_{3,2}}
\def\Hbar{\bar{\H}}
\def\Ch{Chow_{3,2}}

\def\FMps{\overline{{\bf \mcl M}}^{ps}_g}
\def\FMg{\overline{{\bf \mcl M}_g}}

\def\ds{\oplus}

\def\inj{\hookrightarrow}
\def\GL{GL}

\def\Pic{Pic}
\def\PGL{PGL}
\def\cycle{\varpi}

\def\Gr{Gr}

\def\M{\bar{M}}
\def\cO{\mathcal O}
\def\Mps{\overline{M}_3^{ps}}
\def\dps{\delta^{ps}}

\def\Hbar{\H\mod\SL(6)}

\def\cM{\mathcal{M}}
\def\cMps{\bar{\cM}_3^{ps}}
\def\nef{Nef}
\def\lhs{\la^{hs}}
\def\dhs{\d^{hs}}
\def\Mhs{\overline{M}_3^{hs}}
\def\Mcs{\overline{M}_3^{cs}}
\def\HP{\bar{\mcl P}_4}

\def\Gr{Gr}
\def\bar{\overline}

\input xy
\xyoption{all}

\input epsf
\begin{document}

\title[Log minimal model program for $\M_3$]
{Log minimal model program for the moduli space of stable curves
of genus three}

\author[D. Hyeon]{Donghoon Hyeon}
\author[Y. Lee]{Yongnam Lee}

\address[DH]{Department of Mathematical Sciences, Northern Illinois
University, DeKalb, IL 60115, U.S.A., email: hyeon@math.niu.edu}
\address[YL]{ Department of Mathematics, Sogang Univeristy,
Sinsu-dong, Mapo-gu, Seoul 121-742, Korea, email:
ynlee@sogang.ac.kr}

\begin{abstract} In this paper, we completely work out the log
minimal model program for the moduli space of stable curves of
genus three. We employ a rational multiple $\alpha\delta$ of the
divisor $\delta$ of singular curves as the boundary divisor,
construct the log canonical model for the pair $(\bar{\mathcal
M}_3, \alpha\delta)$ using geometric invariant theory as we vary
$\alpha$ from one to zero, and give a modular interpretation of
each log canonical model and the birational maps between them. By
using the modular description, we are able to identify all but one
log canonical models with existing compactifications of $M_3$,
some new and others classical, while the exception gives a new
modular compactification of $M_3$.
\end{abstract}

\maketitle

\tableofcontents

\section{Introduction and Preliminaries}
In \cite{Has}, Hassett gave an outline of a general program
whose goal is to describe the canonical model
\[
\Proj \ds_{m\ge 0} \Gamma(\Mg, m K_{\Mg})
\]
of the moduli space $\Mg$ of stable curves of genus $g$. The main
idea of tackling the problem is to interpolate with the log
canonical models
\[
\Mg(\a) := \Proj \ds_{m\ge 0} \Gamma(\Mg, m(K_{\FMg}+\a \d))
\]
where $\d$ is the divisor on the moduli stack of the singular
curves. We decrease $\a$ from $1$ to $0$ and try to describe the
log canonical models and the relation between the models.

In this article, we carry out the program for $\M_3$: We find all
log canonical models and give a modular interpretation by
realizing them as GIT quotients of suitable Hilbert scheme or Chow
variety of curves. By using the modular interpretation, we are
also able to identify all but one of the log canonical models with
known moduli spaces. The sole exception is, to our knowledge, a
new compactification of $M_3$.

By a theorem of Cornalba and Harris, $K_{\FMg} + \a \d$ is
ample for $9/11 < \a \le 1$ and $\Mg(\a)$ is isomorphic to $\Mg$
for $\a$ in that range. At $\a = 9/11$, it is shown in \cite{HH1}
that the locus of elliptic tails
gets contracted, resulting the moduli space $\Mg^{ps}$ of pseudo-stable
curves of Schubert.

\begin{figure}[!ht]
\centerline{\scalebox{0.4}{\psfig{figure=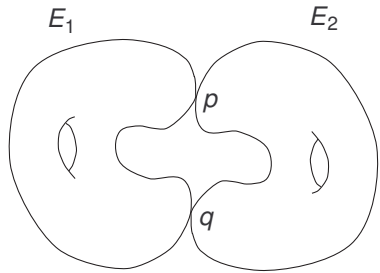}}} \caption{A
generic elliptic bridge}\label{F:eb}
\end{figure}

Our first main theorem is:

\begin{thm}\label{T:main1}
\begin{enumerate}
\item There is a small contraction $\Psi : \Mps \to \M_3(7/10)$ contracting
the locus of {\it elliptic bridges}, and $\M_3(7/10)$ is
isomorphic to the GIT quotient $\Ch\mod\SL(6)$ of the Chow variety
of bicanonical curves;

\item There exists a flip $\Psi^+ : (\Mps)^+ \to \M_3(7/10)$, and $(\Mps)^+$ is
isomorphic to $\M_3(\a)$ for $17/28 < \a < 7/10$. Moreover, $\M_3(\a)$ for
$\a$ in that range is isomorphic to the GIT quotient $\H\mod\SL(6)$ of the Hilbert
scheme of bicanonical curves;

\item There is a divisorial contraction $\Theta : \Mhs \to  \M_3(17/28)$ that contracts
the hyperelliptic locus, and the log canonical model is isomorphic to the compact
 moduli space $\bar Q:= \P(\G(\cO_{\P^2}(4)))\mod \SL(3)$ of plane quartics.
 \end{enumerate}
 \end{thm}
\noindent Here, an {\it elliptic bridge} of genus three is a curve consisting of
two elliptic curves meeting each other in two nodes
(Figure~\ref{F:eb}). Abusing terminology, an elliptic bridge shall sometimes mean
an elliptic component meeting the rest of the curve in two nodes.

\medskip

The second main theorem is the description of the GIT quotient
spaces as moduli spaces. We first introduce two new stability
notions for curves: A complete curve $C$ is said to be {\it
c-stable} if
\begin{enumerate}
\item $C$ has nodes, cusps and tacnodes as singularities;
\item $\o_C$ is ample;
\item a genus one subcurve meets the rest of the
curve in at least two points.
\end{enumerate}
A \emph{cusp} (resp. {\it tacnode}) is the singularity which is locally
analytically isomorphic to $y^2 = x^3$ (resp. $y^2 = x^4$). If a
c-stable curve is not an elliptic bridge, we call it {\it h-stable}.

\begin{thm}\label{T:main2} Let $C$ be a bicanonical curve of genus three.

\begin{enumerate} \item $C$
is Hilbert semistable   if and only if it
is h-stable. It is Hilbert stable  if and only if it
is h-stable and has no tacnode.

\item $C$ is  Chow semistable  if and only if it
is c-stable. It is Chow stable if and only if it is c-stable and
has no tacnode or elliptic bridge. Moreover, all Chow strictly semistable curves
are identified in the moduli space of c-stable curves.
\end{enumerate}
\end{thm}
Due to the theorem, the GIT quotient space $\H\mod\SL(6)$ (resp. $\Ch\mod\SL(6)$)
is the moduli spaces of h-stable curves (resp. c-stable curves),
and we shall denote it by $\M_3^{hs}$ (resp. $\M_3^{cs}$).
The semistability statement of this theorem is proved in \cite{HH2} for
$g \ge 4$. Although some arguments used in the higher genera case go
through, the GIT is quite different in the genus three case. For
instance, for $g \ge 4$, there are many  h-stable curves with tacnodes that are
Hilbert stable, including all irreducible ones.
 Also, the fact that all tacnodal
curves are identified in  $\Mcs$
 allows us to identify
$\Mcs$ with the compact moduli space of plane quartics constructed
by Kondo using the period domains of K3 surfaces
(Proposition~\ref{P:Kondo}).

The following diagram gives a hawk's eye view of the  main
theorems: \small
\[
\xymatrix{ \M_3(1)\simeq\M_3 \ar[d]_-{T} &&&  \\
\M_3(\frac9{11})\simeq  \M_3^{ps} \ar[ddr]_-{\Psi} &
& \M_3(\frac7{10}-\e) \simeq \M_3^{hs}\ar[ddl]^-{\Psi^+}\ar[ddr]_-{\Theta} & & \\
\\
& \M_3(7/10)\simeq \M_3^{cs} & & \M_3(\frac{17}{28})\simeq\bar Q\ar[d]  \\
&&& \star
 }
\]
\normalsize

This program was initiated by Brendan Hassett and Sean Keel, and
the ideas were further developed in  \cite{HH1}. The genus two
case was completely worked out by Hassett in \cite{Has} (see also
\cite{HL1}) and the first couple of steps of the program for the
higher genera case were completed in \cite{HH1} and \cite{HH2}.

\medskip

We work over an algebraically closed field $k$ of characteristic
zero.

\begin{ack} D.H. would like to thank Brendan Hassett for suggesting this
problem and for many helpful conversations. He was partially
supported by KIAS. He gratefully acknowledges Bumsig Kim for the
hospitality and for useful conversations.

 Y.L. would like to
thank Shigeyuki Kondo for his explanation of his compact moduli
space, Shigefumi Mori for his explanation of birational geometry,
Shigeru Mukai for his valuable comments on the geometric invariant
theory. This article was completed while Y.L. was visiting RIMS at
Kyoto University by the JSPS Invitation Fellowship Program, and
Mittag-Leffler Institute. He thanks Shigefumi Mori for the
invitation and hospitality during his stay at RIMS, and Carel
Faber for the invitation and hospitality during his stay at
Mittag-Leffler Institute. He was partially supported by CQUeST
(Sogang University) where he holds a joint appointment.

Both authors thank Noboru Nakayama for his valuable comments and
suggestions that improved this article. Part (3) of Theorem~1 and
the use of Moriwaki divisor in Proposition~\ref{P:morphism} were
suggested to us by Brendan Hassett, and we thank him for generously
sharing his ideas. After our work was completed, we learned that
Sean Keel had obtained some of the results in this paper
independently a few years ago. Both authors were supported by the
Korea Research Foundation Grant funded by the Korean Government
(MOEHRD) (KRF-2005-042-C00005).

\end{ack}

\section{GIT of bicanonical curves of genus
three}\label{S:GIT}

Let $V$ be a vector space of dimension 6. Let $\hilb$ denote the
Hilbert scheme parametrizing subschemes of $\P(V)$ with Hilbert
polynomial $P(m) = 8 m - 2$, and let  $\H$ denote the closure in
$\hilb$ of the locus of the bicanonical images
\[ \xymatrix{ C \ar@{^(->}[r]^-{|\o_C^{\ten 2}|} & \P(V)}
\]
of smooth curves of genus three.
 We shall also consider the corresponding Chow variety $\Ch$ of
bicanonical curves.
 Since $\o_C^{\ten 2}$ is very ample if $C$ is
\emph{c-stable} \cite{HH2}, $\H$ and $\Ch$ include c-stable
curves. We have a natural action of $\SL(V)$ on $\H$ and $\Ch$,
and the aim of this section is to construct the GIT quotients
$\H\mod\SL(V)$ and $\Ch\mod\SL(V)$. Of course we need to specify
the line bundles that we use to linearize the group action: The
Chow variety has a canonical polarization, as it is canonically a
closed subscheme of $\displaystyle{\P(\bigotimes^2 Sym^{8} V^*)}$;
On the Hilbert scheme, we use the line bundle that comes from the
embedding
\[
\phi_m: \H \inj \Gr(P(m), Sym^m V) \inj
\P(\bigwedge^{P(m)}Sym^mV), \quad m \gg 0.
\]
We shall denote the Hilbert point in $\hilb$ of $C$ by $[C]$; Its
image $\phi_m([C])$ will be denoted by $[C]_m$ and will be called
the \emph{$m$th Hilbert point} of $C$.

\begin{defn} (1) \, $C$ is \emph{$m$-Hilbert (semi)stable} if $[C]_m$
is (semi)stable with respect to the natural $\SL(V)$ action on
$\P(\bigwedge^{P(m)}Sym^mV)$. It is said to be \emph{Hilbert
(semi)stable} if it is $m$-Hilbert (semi)stable for all $m \gg 0$.

\noindent (2) \, $C$ is \emph{Chow (semi)stable} if the Chow point
of $C$ is (semi)stable with respect to the natural $\SL(V)$ action
on  $\P(\bigotimes^2 Sym^{8} V^*)$.
\end{defn}

\subsection{Stability computation via Gr\"obner basis}
In this section, we give an overview of some results from
\cite{HHL}. Proofs and detailed analysis as well as some computer
algebra system implementations can be found in that paper.

Let $X \subset \P^N$ be a projective variety and let $I_X$ and $P$
denote the homogeneous ideal and the Hilbert polynomial of $X$.
Let $\rho: \bG_m \to GL(N+1)$ be a one parameter subgroup defined
by $\rho(\a).x_i = \a^{r_i}x_i$ where $x_i$ are  homogeneous
coordinates. The associated one-parameter subgroup of $SL(N+1)$
with weights $(N+1)r_i - \sum r_i$ will be denote by $\rho'$.

\begin{prop}\label{P:algorithm}
 The Hilbert-Mumford index of the $m$th Hilbert point
of $X$ with respect to $\rho'$ is given by the following formula:
\begin{equation}\label{E:mu}
\mu([X]_m, \rho') = -(N+1) \sum_{i=1}^{P(m)} wt_\rho(x^{a(i)}) +
m\cdot P(m)\cdot \sum_{i=0}^Nr_i
\end{equation}
Here, $x^{a(i)}$ are the degree $m$ monomials not in the initial
ideal of $I_X$ with respect to the $\rho$-weighted graded
lexicographic order.
\end{prop}

From the functoriality of the Hilbert-Mumford index and a careful
analysis of the tautological ring of the Hilbert scheme \cite{KnM}
along with the cohomology vanishing property of c-stable curves
\cite{HH2}, one can deduce:

\begin{prop}\label{P:all-m} Let $C \subset \P(V)$ be a bicanonical c-stable curve and
$C^\star$ denote the curve to which $\rho(\a).C$ specializes. If
$C^\star$ is also a bicanonical c-stable curve, then for all $m\ge
2$, $C$ is

\begin{enumerate}
\item $m$-Hilbert stable if and only if $\mu([C]_3,\rho) \ge 2
\mu([C]_2,\rho) > 0 $;

\item $m$-Hilbert strictly semistable if and only if
$\mu([C]_3,\rho) = \mu([C]_2,\rho) = 0$;

\item $m$-Hilbert unstable if and only if   $\mu([C]_3,\rho) \le
2 \mu([C]_2,\rho) < 0 $.
\end{enumerate}
\end{prop}

\subsection{Unstable curves}
As with most GIT problems, it is much easier to classify the
unstable bicanonical curves. We start with the curves that are
easiest to destabilize:

\begin{prop}\label{P:unstable-Chow}
Let $C$ be a  bicanonical curve of genus $\ge 3$ which is not a
curve consisting of two elliptic curves meeting in one tacnode. It
is Chow unstable if it is not c-stable.
\end{prop}
For non c-stable curves that is not the special curve excluded in
the proposition, it is not very difficult to find a destabilizing
one parameter subgroup. The proof given in \cite{HH2} is indeed
valid for all genus $g \ge 2$.

\medskip

By definition, if a non h-stable curve $C$  is not an elliptic
bridge, then it is not c-stable. By
Proposition~\ref{P:unstable-Chow}, $C$ is Chow unstable if it is not
the special curve excluded in Proposition~\ref{P:unstable-Chow}.
With this observation, we can rule out almost all non h-stable
curves as Hilbert unstable, due to the following:

\begin{prop}\label{P:link}(see, e.g. \cite{HH2})
If a projective variety is Chow stable (resp. unstable) then it is
$m$-Hilbert stable (resp. unstable) for $m \gg 0$.
\end{prop}
The following corollary is immediate:
\begin{coro}\label{P:h-unstable}\cite{HHL} If a bicanonical curve $C$ is not
h-stable and is not an elliptic bridge, then it is Hilbert unstable.
\end{coro}

We are left with two cases: (a) Chow instability of a curve
consisting of two elliptic curves meeting in one tacnode; (b)
Hilbert instability  of elliptic bridges.

\begin{prop}\label{P:2g1}
 A bicanonical curve $C$ consisting of two elliptic curves
$E_1$ and $E_2$ meeting in one tacnode is Chow unstable.
\end{prop}

\begin{proof} Suppose that $C$ is Chow semistable.
 The Deligne-Mumford stabilization $D$ of $C$
 consists of an elliptic curve
bridging $E_1$ and $E_2$. But $D$ is also the stabilization of an
irreducible curve $C'$ with two cusps which is Hilbert stable
(\cite{Sch}, see also the proof of Proposition~\ref{P:h-stable}) and
hence Chow stable by Proposition~\ref{P:link}. Since $C'$ and $C$
are both semistable and have the same Deligne-Mumford stabilization,
they must be identified in the quotient space $\Ch\mod SL(V)$. This
contradicts that $C'$ is Chow stable.
\end{proof}

The instability of elliptic bridges is rather nontrivial and we
devote the rest of the section to its proof.

\begin{defn} A {\it snowman} is a genus one curve consisting
of two smooth rational curves meeting in one tacnode.
\end{defn}

\begin{figure}[h]
\centerline{\scalebox{0.4}{\psfig{figure=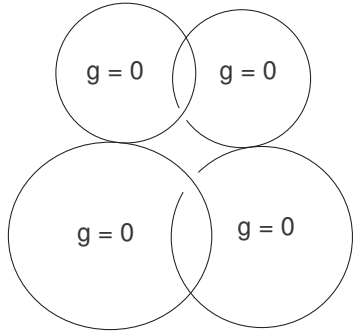}}}
\caption{The conjoined snowmen $C^\ast$}\label{F:snowmen}
\end{figure}

\begin{prop}\label{P:eb-unstable} Let $C^\ast$ be the bicanonical
 curve consisting
of two snowmen meeting in two nodes such that each rational
component meets the rest of the curve in precisely one node and
one tacnode (Figure~\ref{F:snowmen}). Then for every bicanonical
elliptic bridge $C$, there is a one-parameter subgroup $\rho :
\bG_m \to \SL(V)$ such that
\begin{enumerate}
\item  $\rho(\a).C$ specializes to $C^\ast$; \item $C^\ast$ is
Hilbert unstable with respect to $\rho$.
\end{enumerate}
\end{prop}

\begin{proof} Let $C$ be an elliptic bridge consisting of
elliptic curves $E_1$ and $E_2$ meeting in two nodes $p$ and $q$.
Restricting the dualizing sheaf, we have $\o_C|_{E_i} = \o_{E_i}(p +
q)$. This means that each elliptic component  is embedded in $\P^3$
by the $|2 p + 2 q|$ linear system. A generic elliptic curve
embedded by $|2p + 2q|$ linear system is a complete intersection cut
out by
\begin{equation}\label{E:2p2q}
\begin{cases}
x_0x_1 - x_2x_3 = 0 \\
x_1^2 + x_0x_3 + x_0x_2 + a \, x_0x_1 + b \, x_0^2 = 0
\end{cases}
\end{equation}
for some $a, b \in k$. Its $j$-invariant can be obtained by
pulling back the equation (\ref{E:2p2q}) to $\P^1 \ex \P^1$ and
realizing it as a double cover of the $[x,y]$-line ramified over
$\{(x^2 + a x y)^2 - 4 y^2(x y + b x^2) = 0\}$ \footnote{We
learned this computation from N. Nakayama.}:
\[
j = \frac{- 2^83^3 \left(a^2 - 12b\right)^3}{4\left(a^2 -
12b\right)^3 + 27 \left(2a^3 - 72 a b - 4^2 3^3\right)}.
\]

$C \subset \P^5$ comprises of two elliptic curves embedded in such a
fashion:
\[
E_1 : \quad \begin{cases} x_4 = x_5 = 0\\
x_0x_1 - x_2x_3 = 0 \\
x_1^2 + x_0x_3 + x_0x_2 + a \, x_0x_1 + b \, x_0^2 = 0
\end{cases}
\]
\[
E_2 : \quad \begin{cases} x_0 = x_1 = 0\\
x_5x_4 - x_3x_2 = 0 \\
x_4^2 + x_5x_2 + x_5x_3 + a' \, x_5x_4 + b' \, x_5^2 = 0
\end{cases}
\]
Let $\rho : \bb G_m \to \SL(6)$ be the one parameter subgroup
with weights $\{0,1,2,2,1,0\}$. We find that the flat limit
$C^\flat$ of the family $\{\rho(\a).C \}$ at $\a = 0$ does not
depend on the $j$-invariants of $E_1$ and $E_2$,  and is cut out
by the ideal
\begin{equation}\label{E:2snowmen}
I_{C^\flat} = \langle x_1x_5, \,  x_0x_5, \,  x_4^2 + x_2x_5 +
x_4x_5, \,  x_1x_4, \,  x_0x_4, \,  x_2x_3, \,  x_1^2 + x_0x_2 +
x_0x_3\rangle
\end{equation}
which is reduced and has the following associated primes
\[
\left\{ \begin{array}{c} \langle x_{2}, \,x_{1}, \,x_{0},
\,x_{4}^{2}+x_{3} x_{5}\rangle, \, \langle x_{3}, \,x_{1}, \,x_{0},
\,x_{4}^{2}+x_{2} x_{5}\rangle, \\
 \langle x_{5}, \,x_{4}, \,x_{3},
\,x_{1}^{2}+x_{0} x_{2}\rangle, \, \langle x_{5}, \,x_{4}, \,x_{2},
\,x_{1}^{2}+x_{0} x_{3}\rangle \end{array}\right\}.\]
 Note that these are rational conics and each
component meets the rest of the curve in precisely one tacnode and
one node. Hence (\ref{E:2snowmen}) is  the conjoined snowmen curve
$C^\ast$.

We compute the Hilbert-Mumford index of $C^\ast$ with respect to
$\rho$ by using Proposition~\ref{P:algorithm}. For $m = 2$, the
degree two monomials not in $I_{C^\ast}$ are:
\[
x_0^2, \, x_0x_2, \, x_0x_3, \, x_1^2, \, x_1x_2, \, x_1x_3, \,
x_1x_4, \, x_1x_5, \, x_2^2, \, x_2x_4, \, x_3^2, \, x_3x_5, \,
x_4^2, \, x_5^2
\]
of which $\rho$-weights sum up to 30. We have
\[
\mu([C^\ast]_2, \rho) = - 6\cdot 30 + 2\cdot P(2) \cdot (2+2+1+1)
= -180 + 168 = -12.
\]
For $m = 3$, the degree three monomials not in $ I_{C^\ast}$ are
\[
\begin{array}{c}
x_0^3, \,x_0^2x_2, \,x_0^2x_3, \,x_0x_2^2, \,x_0x_3^2, \,x_1^3
x_1^2x_2, \,x_1^2x_3, \,x_1x_2^2, \,x_1x_2x_4, \,x_1x_3^2, \\
\,x_1x_3x_5, \,x_1x_4^2, \,x_1x_5^2, \,x_2^3, \,x_2^2x_4,
\,x_2x_4^2, \,  x_3^3, \,x_3^2x_5, \,x_3x_5^2, \,x_4^3, \,x_5^3
\end{array}
\]
of which $\rho$-weights sum up to 70. The Hilbert-Mumford index
for the third Hilbert point is
\[
\mu([C^\ast]_3, \rho) = - 6\cdot 70 + 2\cdot P(3) \cdot (2+2+1+1)
= -420 + 396 = -24.
\]
Since $C$ and $C^\ast$ are c-stable and $\mu([C^\ast]_3), \rho) =
2 \mu([C^\ast]_2), \rho) < 0$, it follows from
Proposition~\ref{P:all-m} that $[C^\ast]_m$ is unstable with
respect to $\rho$ for all $m \ge 2$.
\end{proof}

\subsection{Semistability proof}
If an h-stable curve $C$ does not have a tacnode, then it is
either Deligne-Mumford stable or pseudo-stable, and the Hilbert
semistability of $C$ essentially follows from \cite{Mum},
\cite{Sch}.

 To describe our strategy of treating the tacnodal ones,
 we first introduce the following two curves which are special
 among all h-stable curves:
\begin{defn}
\begin{enumerate}
\item A \emph{cat-eye} is a genus three curve consisting of two
smooth rational curves $C_1$ and $C_2$ meeting in two tacnodes;
\item An {\it ox} is a genus three curve consisting of three
smooth rational curves $C_1, C_2, C_3$ such that $C_1$ and $C_2$
meet each other in a node and meet $C_3$ in  a tacnode.
\end{enumerate}
\end{defn}

\begin{defn} Let $X \subset \P^N$ be a projective variety and
$\rho : G \to GL(N+1)$ be a one-parameter subgroup. The
\emph{basin of attraction} of a $\rho$-fixed point $x^\star \in X$
is the set
\[
A_\rho(x^\star) := \{ x \in X \, | \, \rho(\a).x \rightsquigarrow
x^\star\}.
\]
\end{defn}
It is easy to show that if $x^\star$ is strictly semistable with
respect to $\rho$, then $x^\star$ is semistable if and only if every
$x \in A_\rho(x^\star)$ is semistable.  We shall enumerate all
possible h-stable replacements of a cat-eye (resp. the ox) by
deformation theory, and show that those are in a basin of attraction
of the particular cat-eye (resp. the ox) by direct flat limit
computation.

\begin{prop} A genus three h-stable curve $C$ has infinite automorphisms
if and only if it is a {\it cat-eye}  or an {\it ox}
(Figure~\ref{F:h-stable-tac}).
\end{prop}

\begin{proof} An h-stable curve without a
tacnode is Deligne-Mumford stable or pseudo-stable, and has finite
automorphisms. Since an automorphism of $C$ lifts to an
automorphism of its normalization, it is easy to see that among
the tacnodal h-stable curves  listed in
Figure~\ref{F:h-stable-tac}, only (d) and (i) have infinite
automorphisms.

%\begin{figure}[htb]
%\centerline{\scalebox{0.7}{\psfig{figure=semistable-genus3-new-1.eps}}}
%\caption{h-stable curves of genus three with
%tacnodes}\label{F:h-stable-tac}
%\end{figure}

\begin{figure}[!t]
  \begin{center}
    \includegraphics[width=6.5in]{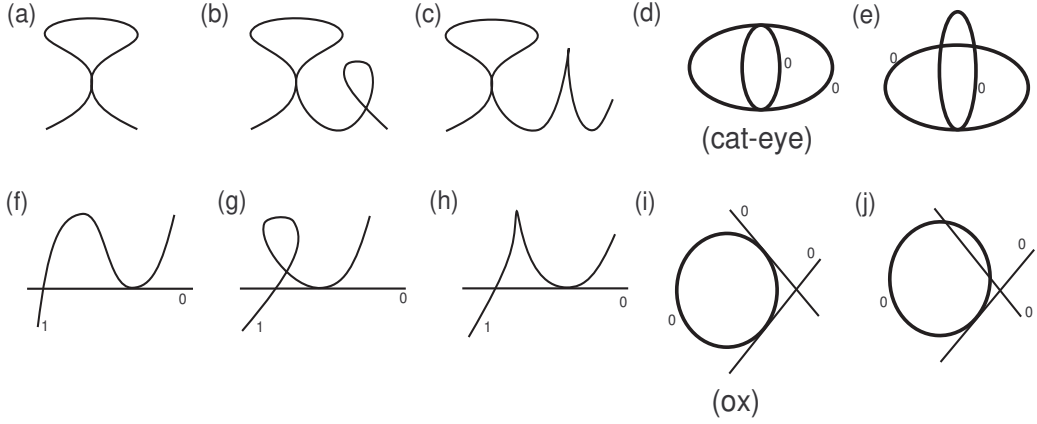}
  \end{center}

  \caption{h-stable curves of genus three with tacnodes}
  \label{F:h-stable-tac}
\end{figure}
\end{proof}

Let $C^{at}$ be a cat-eye consisting of two smooth rational curves
$C_1$ and $C_2$ meeting in two tacnodes $p$ and $q$, embedded in
$\P^5$ by the bicanonical system $|\o_{C^{at}}^{\ten 2}|$. There
is a one parameter family $\{C^{at}_\beta\}$ of cat-eyes where the
parameter $\beta$ encodes the tangent space identifications. Since
$\o_{C^{at}}|_{C_i} \simeq \o_{C_i}(2p + 2q)$, each $C_i$ is the
second Veronese image of a plane conic. Hence a cat-eye is the
second Veronese image of the plane quartic
\begin{equation}\label{E:plane-cat-eye}
\{ (x_1^2 + x_0x_2)(\beta x_1^2 + x_0x_2) = 0 \}, \quad \beta \in
k\setminus\{0, 1\}.
\end{equation}
which has tacnodes at $[1,0,0]$ and $[0,0,1]$. Note that $C_1$
(resp. $C_2$)  has a parametrization coming from $[s, t] \mapsto
[s^2, s t, - t^2]$ (resp., $[s^2,  s t, - \beta t^2]$) and
automorphisms $[s, t] \mapsto [\a s, t]$, $\a \in \bG_m$. This
induces an action of the one-parameter subgroup $ \rho : \bG_m \to
\GL(6)$ defined by $\rho(\a)\cdot x_i = \a^{r_i}x_i$ with the
weights $(r_0, \cdots, r_5) = (4,2,0,3,2,1)$.

\begin{prop}\label{P:cat-ss} A cat-eye is $m$-Hilbert strictly semistable with
respect to $\rho$ for all $m \ge 2$.
\end{prop}

\begin{proof} The Gr\"obner basis of $C$ with respect to the weighted GLex
is
\begin{equation}\label{E:gb-cat-eye}
\begin{array}{c}
x_{1} x_{{2}}-x_{{5}}^{2}, \,
      x_{{2}} x_{{3}}-x_{{4}} x_{{5}}, \,
      x_{1} x_{{4}}-x_{{3}} x_{{5}}, \,
      x_{0} x_{{2}}-x_{{4}}^{2}, \,
      x_{1}^{2} {\beta}+x_{{3}} x_{{5}} {\beta}+x_{{3}}
      x_{{5}}+x_{{4}}^{2}, \\
      x_{1} x_{{5}}^{2} {\beta}+x_{{4}} x_{{5}}^{2} {\beta}+x_{{2}} x_{{4}}^{2}+x_{{4}} x_{{5}}^{2}, \,
      x_{{2}} x_{{4}} x_{{5}}^{2} {\beta}+x_{{5}}^{4} {\beta}+x_{{2}}^{2} x_{{4}}^{2}+x_{{2}} x_{{4}} x_{{5}}^{2},
      \, x_{0} x_{{5}}-x_{{3}} x_{{4}}, \,
       \\
      x_{1} x_{{3}} x_{{5}} {\beta}+x_{{3}} x_{{4}} x_{{5}} {\beta}+x_{{3}} x_{{4}} x_{{5}}+x_{{4}}^{3}, \,
      x_{{3}} x_{{5}}^{3} {\beta}+x_{{4}}^{2} x_{{5}}^{2} {\beta}+x_{{2}} x_{{4}}^{3}+x_{{4}}^{2} x_{{5}}^{2}, \,
      x_{0} x_{1}-x_{{3}}^{2},\\
      x_{{3}}^{2} x_{{5}}^{2} {\beta}+x_{{3}} x_{{4}}^{2} x_{{5}} {\beta}+x_{{3}} x_{{4}}^{2} x_{{5}}+x_{{4}}^{4}, \,
      x_{{3}}^{3} x_{{5}} {\beta}+x_{{3}}^{2} x_{{4}}^{2} {\beta}+x_{0} x_{{4}}^{3}+x_{{3}}^{2} x_{{4}}^{2}, \, \\
      x_{0} x_{{3}}^{2} x_{{4}} {\beta}+x_{{3}}^{4} {\beta}+x_{0}^{2} x_{{4}}^{2}+x_{0} x_{{3}}^{2} x_{{4}}, \,
       x_{1} x_{{3}}^{2} {\beta}+x_{{3}}^{2} x_{{4}} {\beta}+x_{0} x_{{4}}^{2}+x_{{3}}^{2}x_{{4}}.\\
 \end{array}
\end{equation}
Hence the degree two monomials not in the initial ideal are
\begin{equation}\label{E:inI-2-cat-eye}
x_0^2, \,
      x_0 x_3, \,
      x_0 x_4, \,
      x_1 x_3, \,
      x_1 x_5, \,
      x_2^2, \,
      x_2 x_4, \,
      x_2 x_5, \,
      x_3^2, \,
      x_3 x_4, \,
      x_3 x_5, \,
      x_4^2, \,
      x_4 x_5, \,
      x_5^2
\end{equation}
and the sum of the weights of these monomials is $56$. On the other
hand, the average weight (the second term of the right hand side of
(\ref{E:mu}) divided by $N+1 = 6$)  is $\frac{2 \cdot 14\cdot 12}{6}
= 56$, and we have $\mu([C]_2, \rho) = 0$ by
Proposition~\ref{P:algorithm}.

The degree three monomials in the initial ideal are multiples of
(\ref{E:inI-2-cat-eye}) together with $ x_1x_5^2$, $ x_1x_3x_5$,
and $x_1x_3^2$. Hence the degree three monomials not in the
initial ideal are
\begin{equation}\label{E:inI-3-cat-eye}
\begin{array}{c}
x_{0}^{3}, \,   x_{0}^{2} x_{3}, \,   x_{0}^{2} x_{4}, \,   x_{0}
x_{3}^{2}, \, x_{0} x_{3} x_{4}, \,   x_{0} x_{4}^{2}, \,
x_{2}^{3}, \, x_{2}^{2} x_{4}, \,
 x_{2}^{2} x_{5}, \,    x_{2} x_{4}^{2}, \,    x_{2} x_{4} x_{5}, \, \\
  x_{2} x_{5}^{2}, \,  x_{3}^{3}, \,  x_{3}^{2} x_{4}, \,  x_{3}^{2} x_{5}, \,
  x_{3} x_{4}^{2}, \,      x_{3} x_{4} x_{5}, \,      x_{3} x_{5}^{2}, \,
    x_{4}^{3}, \,      x_{4}^{2} x_{5}, \,      x_{4} x_{5}^{2}, \,
     x_{5}^{3}\\
\end{array}
\end{equation}
of which weights sum up to $132$. The average weight is $\frac{3
\cdot 22\cdot 12}{6} = 132$, and $C$ is 3-Hilbert strictly
semistable. The assertion now follows from
Proposition~\ref{P:all-m}.
\end{proof}

Consider the irreducible h-stable tacnodal curves (a)$\sim$(c).
These are obtained from a genus one curve by identifying two points
and their tangent lines. Given a genus one curve $E$ and two
distinct simple points $q$ and $r$, one can construct a
one-parameter family $\{E_{\beta}\}$ of tacnodal curves by
identifying $q$, $r$, and identifying the tangent lines by
\[
\frac{\partial}{\partial \sm_q} = \beta \frac{\partial}{\partial
\sm_r}
\]
where $\sm_q$ and $\sm_r$ denote the local parameters at $q$ and
$r$.

\begin{prop}\label{P:cat-flat} The flat limit of
$\{\rho(\a).E_\beta\}$ at $\beta = 0$ is $C^{at}_{\beta}$.
\end{prop}

\begin{proof}
We start with $E \subset \P^3$ embedded by $|2q + 2r|$. This is
defined by
\[
x_0x_1 - x_2x_3 = x_1^2 + x_0x_3 + x_0x_2 + a x_0x_1 + b x_0^2 = 0
\]
where $a, b \in k$. Here $q = [0,0,1,0]$ and $r = [0,0,0,1]$. The
following projection identifies $q$ and $r$:
\begin{equation}\label{E:pr}
pr : [x_0, \cdots, x_3] \mapsto [x_0,x_1,x_2+\beta x_3]
\end{equation}
The order of vanishing of $x_0, \cdots, x_3$ at $q$ and $r$ are
\begin{equation}\label{E:van-o}
\begin{tabular}{|c|c|c|c|c|}
  \hline
   & $x_0$ & $x_1$ & $x_2$  & $x_3$ \\
  \hline
  $ord_q$ & 2 & 1 & 0 & 3 \\
  \hline
  $ord_r$ & 2 & 1 & 3 & 0 \\
  \hline
\end{tabular}
\end{equation}
From this, we deduce that under $pr$, the tangent lines are
identified by the relation
 $\frac{\partial}{\partial \sigma_q}
= \beta \frac{\partial}{\partial \sigma_r}$. The image of $E$ under
$pr$ is the quartic curve defined by:
\small\begin{equation}\label{E:prE}
\begin{array}{c}
f = x_{0}^{2} x_{1}^{2} a^{2} \beta+2 x_{0}^{3} x_{1} {a} {b}
\beta+x_{0}^{4} b^{2} \beta+2 x_{0} x_{1}^{3} {a} \beta+x_{0}^{2}
x_{1} x_{{2}} {a} \beta+2 x_{0}^{2} x_{1}^{2} {b} \beta+x_{0}^{3}
x_{{2}} {b} \beta
\\ +x_{0}^{3} x_{1} \beta^{2} +x_{0}^{2} x_{1} x_{{2}} {a}+x_{0}^{3} x_{{2}} {b}-2 x_{0}^{3} x_{1} \beta+x_{1}^{4} \beta+x_{0} x_{1}^{2} x_{{2}} \beta+x_{0}^{3} x_{1}+x_{0} x_{1}^{2} x_{{2}}+x_{0}^{2} x_{{2}}^{2} = 0.
\end{array}
\end{equation}\normalsize
The second Veronese image of (\ref{E:prE}) is the bicanonical
image of $E_\beta$ of which ideal can be readily computed. We omit
it as we do not need it.

% and its ideal is
%\begin{equation}\label{E:ideal-irred-tac}
%\begin{array}{c}
%x_3 x_4-x_0 x_5,\,x_1 x_4-x_3 x_5,\,x_2 x_3-x_4 x_5,\,x_1 x_2-x_5^2,\,x_0 x_2-x_4^2,\, \\
%x_0 x_1-x_3^2,\,
%x_3^2 a^2 {t}+2 x_0 x_3 {a} {b} {t}+x_0^2 b^2 {t}+2 x_1 x_3 {a} {t}+x_0 x_5 {a} {t}+2 x_3^2 {b} {t} \\
%+x_0 x_4 {b} {t}+x_0 x_3 t^2+x_0 x_5 {a}+x_0 x_4 {b}+x_1^2 {t}-2 x_0 x_3 {t}+x_3 x_5 {t}+x_0 x_3+x_4^2+x_3 x_5,\, \\
%x_0 x_5^2 a^2 {t}+2 x_0 x_4 x_5 {a} {b} {t}+x_0 x_4^2 b^2 {t}+x_4^2 x_5 {a} {t}+2 x_3 x_5^2 {a} {t}+x_4^3 {b} {t} \\
%+2 x_0 x_5^2 {b} {t}+x_0 x_4 x_5 t^2+x_4^2 x_5 {a}+x_4^3 {b}-2 x_0 x_4 x_5 {t}+x_1 x_5^2 {t}+x_4 x_5^2 {t}+x_2 x_4^2 +x_0 x_4 x_5+x_4 x_5^2,\, \\
%x_0 x_3 x_5 a^2 {t}+2 x_0^2 x_5 {a} {b} {t}+x_0^2 x_4 b^2 {t}+2 x_3^2 x_5 {a} {t}+x_0 x_4 x_5 {a} {t}+x_0 x_4^2 {b} {t}\\
%+2 x_0 x_3 x_5 {b} {t}+x_0^2 x_5 t^2+x_0 x_4 x_5 {a}+x_0 x_4^2 {b}-2 x_0^2 x_5 {t}+x_1 x_3 x_5 {t}+x_0 x_5^2 {t}+x_4^3+x_0^2 x_5+x_0 x_5^2,\, \\
%x_4^2 x_5^2 a^2 {t}+2 x_4^3 x_5 {a} {b} {t}+x_4^4 b^2 {t}+x_2 x_4^2 x_5 {a} {t}+2 x_4 x_5^3 {a} {t}+x_2 x_4^3 {b} {t}+2 x_4^2 x_5^2 {b} {t} \\
%+x_4^3 x_5 t^2+x_2 x_4^2 x_5 {a}+x_2 x_4^3 {b}-2 x_4^3 x_5 {t}+x_2 x_4 x_5^2 {t}+x_5^4 {t}+x_2^2 x_4^2+x_4^3 x_5+x_2 x_4 x_5^2
%\end{array}\end{equation}

We shall now prove  that the irreducible tacnodal curve $C$
degenerates to a cat-eye along the action of the one parameter
subgroup with weights $\{0,2,4,1,2,3\}$.\footnote{This appears
slightly different from the weights of $\rho$ in
Proposition~\ref{P:cat-ss}, but this does not affect our argument as
they are projectively equivalent.}  Since taking flat limits
commutes with taking Veronese images, we shall work with the plane
quartic (\ref{E:prE}) and show that it degenerates to the cat-eye
(\ref{E:plane-cat-eye}) along the action of the one parameter
subgroup  with weights $\{0,1,2\}$. By slightly abusing notations,
we let $C$ denote the quartic curve (\ref{E:prE}) and $\rho$, the
one-parameter subgroup of $\GL(3)$ with weights $\{0,1,2\}$. Then
the flat projective closure of $\{\rho(\a). C\}$ is defined by
\[
\begin{array}{l}
\a^4 f(x_0, \a\inv x_1, \a^{-2} x_2) = \a^{2} x_{0}^{2} x_{1}^{2}
a^{2} \beta + \a^{3} 2 x_{0}^{3} x_{1} {a} {b} \beta
+ \a^4 x_{0}^{4} b^{2} \beta+ \a 2 x_{0} x_{1}^{3} {a} \beta \\
+ \a x_{0}^{2} x_{1} x_{{2}} {a} \beta + \a^{2} 2 x_{0}^{2}
x_{1}^{2} {b} \beta + \a^{2} x_{0}^{3} x_{{2}} {b} \beta  + \a^{3}
x_{0}^{3} x_{1} \beta^{2}
+ \a x_{0}^{2} x_{1} x_{{2}} {a}+ \a^{2} x_{0}^{3} x_{{2}} {b} \\
- \a^{3} 2 x_{0}^{3} x_{1} \beta +  x_{1}^{4} \beta+  x_{0}
x_{1}^{2} x_{{2}} \beta + \a^{3} x_{0}^{3} x_{1}+  x_{0} x_{1}^{2}
x_{{2}} + x_{0}^{2} x_{{2}}^{2} = 0,
\end{array}
\]
\normalsize and the flat limit at $\a = 0$ is defined by
\[
x_{1}^{4} \beta+  x_{0} x_{1}^{2} x_{{2}} \beta +   x_{0}
x_{1}^{2} x_{{2}} + x_{0}^{2} x_{{2}}^{2} = (x_1^2 +
x_0x_2)(x_1^2\beta+x_0x_2) = 0
\]
which is precisely the plane cat-eye (\ref{E:plane-cat-eye}).
\end{proof}
There is another curve that specializes to a cat-eye under the
$\rho$-action: A curve of the form (e) is the second Veronese
image of the quartic curve
\[
D := \{ (x_1^2 + x_0x_2) (\b x_1^2 +x_0x_2 + \g x_2^2) = 0 \}, \quad
\g \in \bG_m,
\]
where the parameter $\g$ encodes the cross ratio of the three
points of attachments. Since the $\rho$-weight of the extra term
$\g x_2^2$ is zero, the flat limit of $\{\rho(\a). D\}$ at $\a =
0$ is $\{ (x_1^2 + x_0x_2) (\b x_1^2 +x_0x_2) = 0\}$.

 From Propositions \ref{P:cat-ss} and \ref{P:cat-flat}, it
follows that

\begin{coro}\label{C:abce} A curve of the form (a), (b), (c), (e) is Hilbert semistable
if and only if the corresponding cat-eye is semistable.
\end{coro}

%\begin{figure}[htb]
%\centerline{\scalebox{0.4}{\psfig{figure=tac-degen.eps}}}
%\caption{Degeneration to the cat-eye and the ox}\label{F:tac-degen}
%\end{figure}

\begin{figure}[!t]
  \begin{center}
    \includegraphics[width=5in]{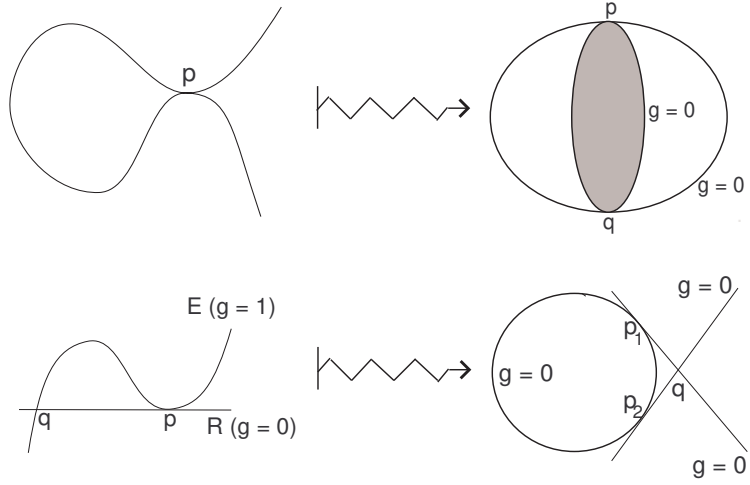}
  \end{center}

  \caption{Degeneration to the cat-eye and the ox}
  \label{F:tac-degen}
\end{figure}

Recall that an ox is a genus three curve consisting of three
smooth rational curves $C_1, C_2, C_3$ such that $C_1$ and $C_2$
meets each other in a node $q$ and $C_3$, in tacnodes $p_1$ and
$p_2$ respectively. Let $C^{ox}$ be an ox, and consider the
restriction of $\o_{C^{ox}}^{\ten 2}$ to each component:
\[
\begin{array}{l}
\o_{C^{ox}}^{\ten 2}|_{C_i} \simeq \o_{C_i}^{\ten 2}(4 p_i + 2 q), \quad i = 1, 2 \\
\o_{C^{ox}}^{\ten 2}|_{C_3} \simeq \o_{C_3}^{\ten 2}(4 p_1+ 4
p_2).
\end{array}
\]
This implies that the bicanonical image of $C^{ox}$ has two smooth
conics $C_1$ and $C_2$ meeting each other in a node and meeting a
smooth rational quartic curve in a tacnode. Since the second
Veronese image $C$ of the plane quartic $\{x_0x_2(x_0x_2 - x_1^2)
= 0 \}$ has precisely such components, and since ox is unique up
to isomorphism, it follows that the bicanonical image of an ox is
projectively equivalent to $C$.

 The plane quartic $\{\,x_0x_2(x_0x_2 - x_1^2) = 0 \,\}$
admits automorphisms $[x_0,x_1,x_2] \mapsto [x_0, \a x_1,
\a^2x_2]$, $\a \in \bG_m$. The second Veronese image $C$ has
associated automorphisms $x_i \mapsto \a^{r_i}x_i$ where
$(r_0,\dots,r_5) = (0,2,4,1,2,3)$. Let $\rho$ be the one parameter
subgroup with these weights.

\begin{prop}\label{P:ox-ss}
 $C$ is $m$-Hilbert strictly semistable with respect to
$\rho$ for all $m\ge 2$.
\end{prop}

\begin{proof} The ideal of $C$ is generated by
\[
\begin{array}{c}
x_{4}^{2}-x_{3} x_{5}, \, x_{3} x_{4}-x_{0} x_{5}, \, x_{1}
x_{4}-x_{3} x_{5}, \,   x_{2} x_{3}-x_{4} x_{5}, \,  \\x_{1}
x_{2}-x_{5}^{2}, \, x_{0} x_{2}-x_{3} x_{5}, \, x_{0}
x_{1}-x_{3}^{2}
\end{array}
\]
The degree two monomials not in the initial ideal are
\[
x_{0}^{2}, \,x_{0} x_{3}, \,x_{0} x_{4}, \,x_{1}^{2}, \,x_{1}
x_{3}, \,x_{1} x_{5}, \,x_{2}^{2}, \,x_{2} x_{4}, \,x_{2} x_{5},
\,x_{3}^{2}, \,x_{3} x_{4}, \,x_{4}^{2}, \,x_{4} x_{5},
\,x_{5}^{2}.
\]
The sum of the weights of these monomials is $56$. This is equal
to the average weight $\frac{2\cdot P(2)\cdot (4+2+0+3+2+1)}{6} =
\frac{2\cdot 14 \cdot 12}{6} $. Therefore $\mu([C]_2, \, \rho) =
0$ by Proposition~\ref{P:algorithm}.

Now we analyze the 3rd Hilbert point of $C$. The degree three
monomials not in the initial ideal are
\[
\begin{array}{c}
x_{0}^{3}, \,x_{0}^{2} x_{3}, \,x_{0}^{2} x_{4}, \,x_{0}
x_{3}^{2}, \,x_{0} x_{3} x_{4}, \,x_{1}^{3}, \,x_{1}^{2} x_{3},
\,x_{1}^{2} x_{5}, \,x_{1} x_{3}^{2}, \,x_{1} x_{5}^{2},
\,x_{2}^{3}, \, \\x_{2}^{2} x_{4}, \, x_{2}^{2} x_{5}, \,x_{2}
x_{4} x_{5}, \,x_{2} x_{5}^{2}, \,x_{3}^{3}, \,x_{3}^{2} x_{4},
\,x_{3} x_{4}^{2}, \,x_{4}^{3}, \,x_{4}^{2} x_{5}, \,x_{4}
x_{5}^{2}, \,x_{5}^{3}
\end{array}
\]
The weights of these sum up to $132$, which is equal to the
average weight $\frac{3\cdot (8\cdot 3 - 2) \cdot 12}{6}$. Hence
$\mu([C]_3,\rho) = 0$. The assertion now follows from
Proposition~\ref{P:all-m}.
\end{proof}

Let $F$ be a bicanonical curve of genus three consisting of a
smooth elliptic curve $E$ meeting a smooth rational curve $R$ in a
tacnode $p$ and a node $q$. We have $\o_F|_{E} = \o_E(2p + q)
\simeq \mcl O(2p + q)$ and $\o_F|_R = \o_R(2p+q)$. Hence $F$ is
the second Veronese image of
\[
F' = \{ f =  x_0 \left(x_1^2 x_2 - x_0 (x_0 - x_2) (x_0 - \ell
x_2)\right) = 0 \} \subset \P^2.
\]

\begin{prop}\label{P:ox-flat} The flat limit of $\{\rho(\a). F\}$
is the ox.
\end{prop}
\begin{proof}
It suffices to show that $F'$ degenerates to the canonical model
of ox along the action of the one-parameter subgroup with weights
$(0, 1, 2)$, which we abuse notation and denote by $\rho$. Then
the flat projective closure of $\{\rho(\a). F\}$ is defined by
\[
\a^4 f(\a^{-2}x_0, \a^{-1}x_1, x_2) = \a^4 x_0^4  - \a^2 x_0^3x_2
\ell - \a^2 x_0^3 x_2 + x_0^2x_2^2 \ell - x_0x_1^2x_2.
\]
At $\a = 0$, this gives
\[
\{ x_0x_2 (x_0x_2 \ell - x_1^2) = 0 \} = \{x_0x_2 (x_0x_2  -
(x_1/\sqrt{\ell})^2) = 0\}
\]
which is the ox.
\end{proof}

It is easy to see that a curve of the form (j) also belongs to the
basin of attraction of the ox with respect to $\rho$: Such a curve
is the second Veronese image of the plane quartic
\[
D = \{ x_0x_2 ( x_0x_2 + x_1^2 + \g x_2^2) = 0 \}.
\]
Since the extra term $\g x_0 x_2^3$ has zero $\rho$-weight, the
flat limit of the family $\{ \rho(\a). D \}$ is $ \{ x_0x_2 (
x_0x_2 + x_1^2) = 0 \}$.

In view of Propositions \ref{P:ox-ss} and \ref{P:ox-flat}, we have
\begin{coro}\label{C:fghj}
Tacnodal curves of the form (f), (g), (h), (j) are Hilbert
semistable if and only if the ox is Hilbert semistable.
\end{coro}
It is easy to deduce from the defining equation
(\ref{E:plane-cat-eye}) that the ox is the flat limit of the family
$\{C^{at}_\beta\}$ at zero and infinity.

\medskip

Now we are ready to complete the proof of the first part of Theorem~\ref{T:main2}:

\begin{prop}\label{P:h-stable} An h-stable bicanonical curve is Hilbert semistable.
Furthermore, an h-stable bicanonical curve without a tacnode is
Hilbert stable.
\end{prop}

\begin{proof} Let $C$ be an h-stable curve. By
a theorem of Mumford (\cite{Mum}, 4.15), a nonsingular curve is
Chow stable and hence Hilbert stable by Proposition~\ref{P:link}.
For singular h-stable curves $C$, we apply a standard degeneration
argument. If
 $C$ does not have a tacnode and is not an elliptic bridge, then it
 does not belong to any basin of attraction. Hence $C$ is the only
 possibly semistable replacement (other than its Deligne-Mumford
 stable model $D$ which is unstable unless $D = C$).
 It follows that $C$ is Hilbert stable.

It remains to show that h-stable curves with tacnodes are Hilbert
semistable.  In view of Propositions \ref{P:cat-flat} and
\ref{P:ox-flat}, it suffices to establish the semistability for
the cat-eyes and the ox. We first note that a generic cat-eye must
be semistable: If not, all cat-eyes and the ox would be unstable
since the ox is in the closure of the locus of cat-eyes, and
consequently bicanonical elliptic bridges would not have Hilbert
semistable stabilizations.

Consider the cat-eye $C^{at}_\beta$  and a smoothing $\pi: \mcl C
\to B$ of it. By the semistable replacement theorem, there is a
family $\pi': \mcl C' \to B$ whose generic fibre $\mcl C'_\eta$ is
isomorphic to the generic fibre $\mcl C_\eta$ of $\mcl C$ and the
special fibre $\mcl C'_s$ is Hilbert semistable. Since the cat-eyes
and the ox are the only ones that have elliptic bridges as the
Deligne-Mumford stabilization, $\mcl C'_s$ must be $C^{at}_{\beta'}$
for some $\beta'$. Here $C^{at}_0 = C^{at}_\infty $ means the ox.

By choosing a basis of $\pi_*(\o_{\mcl C/B})$ (resp.
$\pi'_*(\o_{\mcl C'/B})$), we obtain an embedding $\mcl C \inj
\P^2 \ex B$ (resp. $\mcl C' \inj \P^2 \ex B$). These in turn lead
to morphisms $f, f' : B \to \bar Q:=
\P(\Gamma(\SS_{\P_2}(+4)))\mod \PGL(3)$. Since the generic
fibres of $\mcl C'$ and $\mcl C$ are isomorphic, $f(0) = f'(0)$.
But since the cat-eyes and the ox are separated in the moduli
space $\bar Q$ of plane quartic curves (Section~\ref{S:17/28}), it
follows that $\mcl C_s = \mcl C'_s$ and  $C^{at}_\beta$ is
semistable.
\end{proof}

We let $\Mhs$ denote the GIT quotient $\Hbar$ and call it the moduli
space of h-stable curves of genus three. Since tacnodal curves are
identified with suitable cat-eyes and all cat-eyes are separated in $\Mhs$,

\begin{coro}\label{C:Z^+}
 The locus of tacnodal curves in $\Mhs$ is a smooth rational
curve.
\end{coro}

\begin{figure}[!b]
  \begin{center}
    \includegraphics[width=5in]{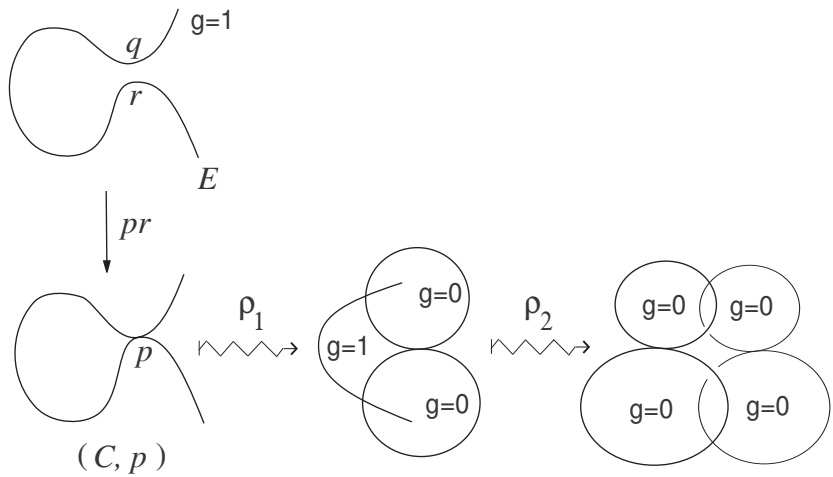}
  \end{center}

  \caption{Degeneration of an elliptic bridge to the conjoined
snowmen}
  \label{F:eb-to-snowmen}
\end{figure}

It remains to show that a c-stable curve is Chow semistable.
Since a c-stable curve that is not an elliptic bridge is h-stable, we
deduce from Propositions \ref{P:link} and \ref{P:h-stable} that
\begin{prop} A bicanonical c-stable curve that is not an elliptic bridge is Chow
semistable. Moreover, such a curve is Chow stable if and only if it has no tacnode.
\end{prop}
To show that elliptic bridges are Chow semistable, we will need
the following  degeneration results.
\begin{prop}\label{P:eb} There are two one-parameter subgroups
$\rho_1$ and $\rho_2$ of $\GL(6)$ such that the following holds
(Figure~\ref{F:eb-to-snowmen}) :

\begin{enumerate}
\item Let $C$ be an irreducible tacnodal curve of genus three.
The flat limit $C^\flat$ of $\{\rho_1(\a). C\}$ at $\a = 0$
consists of an elliptic curve meeting a snowman in two nodes: one
in the head and the other in the body.

\item The flat limit of $\{\rho_2(\a). C^\flat\}$ at $\a = 0$ is
the conjoined snowmen $C^{\ast}$
(Proposition~\ref{P:eb-unstable}).
\end{enumerate}
\end{prop}

%\begin{figure}[h]
%\centerline{\scalebox{0.5}{\psfig{figure=eb-to-snowmen.eps}}}
%\caption{Degeneration of an elliptic bridge to the conjoined
%snowmen}\label{F:eb-to-snowmen}
%\end{figure}

\begin{proof} (a) \,  We retain the homogeneous coordinates $x_0, x_1, x_2$
and equations for $(E,q,r)$ and $E_\beta$ from the proof of
Proposition~\ref{P:cat-flat}. The second Veronese embedding takes
$[x_0, x_1, x_2]$ to $[x_0^2, x_1^2, x_2^2, x_0x_1, x_0x_2,
x_1x_2]$. Let $\rho_1$ be the one-parameter subgroup of $\GL(6)$
with weights $(0,0,2,0,0,1)$. The ideal of the flat limit of
$\{\rho_1(\a). C\}$ at $\a = 0$ is generated by:

\begin{equation}\label{E:Cflat}
\begin{array}{l}
z_3 z_5, \,z_0 z_5, \,z_2 z_3, \,z_1 z_2-z_5^2, \,z_0 z_2, \,z_0
z_1-z_3^2, \, z_1 z_5^2 \beta +z_4 z_5^2 \beta +z_2 z_4^2+z_4
z_5^2,
\, \\
 z_2 z_4
z_5^2 \beta +z_5^4 \beta +z_2^2 z_4^2+z_2 z_4 z_5^2 , \,
 z_1^2 z_5 \beta +z_1 z_4 z_5 \beta +z_1 z_4 z_5+z_4^2 z_5, \, \\

z_3^2 a^2 \beta +2 z_0 z_3 a b \beta +z_0^2 b^2 \beta  +2 z_1 z_3
a \beta  +z_3 z_4 a \beta  +2 z_3^2 b \beta +z_0 z_4 b \beta +z_0
z_3
\beta ^2 \\
+z_3 z_4 a + z_0 z_4 b  +z_1^2 \beta
-2 z_0 z_3 \beta +z_1 z_4 \beta +z_0 z_3+z_1 z_4+z_4^2, \, \\
\end{array}
\end{equation}
This has the following associated primes: \smaller
\begin{enumerate}
\item[$\bullet$] $C_1 := (z_3, \, z_0, z_1+z_4, \, z_2 z_4+z_5^2);$
\item[$\bullet$] $C_2 := (z_3, \,z_0, \,z_1 \beta +z_4, \,z_1 z_2-z_5^2, \,z_5^2 \beta +z_2
z_4); $
\item[$\bullet$] \, $E := (z_5, \,z_2, \,z_0 z_1-z_3^2, \,z_0^2 b^2 \beta +2 z_0 z_3 b
\beta a+z_3^2 \beta  a^2+2 z_3^2 b \beta  +z_0 z_4 b \beta +z_0
z_3 \beta
^2 \\
+2 z_1 z_3 \beta  a+z_3 z_4 \beta a+z_0 z_4 b+z_1^2 \beta -2 z_0
z_3
\beta +z_1 z_4 \beta +z_3 z_4 a+z_0 z_3 + z_1 z_4+z_4^2, \, \\
z_0 z_3^2 b^2 \beta +2 z_3^3 b \beta a+z_1 z_3^2 \beta  a^2 +2 z_1
z_3^2 b \beta +z_3^2 z_4 b \beta +z_3^3 \beta ^2+2 z_1^2 z_3 \beta
a
 +z_1 z_3 z_4 \beta  a  \\
 +z_3^2 z_4 b+z_1^3 \beta -2 z_3^3 \beta
+z_1^2 z_4 \beta +z_1 z_3 z_4 a+z_3^3+z_1^2 z_4+z_1 z_4^2, \,z_3^4
b^2 \beta +2 z_1 z_3^3 b \beta a \\
+z_1^2 z_3^2 \beta  a^2+2 z_1^2 z_3^2 b \beta +z_1 z_3^2 z_4 b
\beta +z_1 z_3^3 \beta ^2 +2 z_1^3 z_3 \beta a  +z_1^2 z_3 z_4
\beta a
+z_1 z_3^2 z_4 b +z_1^4 \beta  \\
-2 z_1 z_3^3 \beta +z_1^3 z_4 \beta +z_1^2 z_3 z_4 a+z_1
z_3^3+z_1^3 z_4+z_1^2 z_4^2).$
\end{enumerate}
\normalsize The first two components are conics and the third is
an elliptic curve isomorphic to $E$. The two conics meet in a
single tacnode at $[0, 0, 1, 0, 0, 0]$ and the elliptic component
intersects $C_1$ and $C_2$ in nodes at $[0,1,0,0,-1,0]$ and
$[0,1,0,0,-\beta, 0]$, respectively.

\

\noindent (2) \,  Let $\rho_2$ be the one-parameter subgroup with
weights $(0,2,2,1,2,2)$. The flat limit of $\{\rho_2(\a).
C^\flat\}$ at $\a = 0$ is defined by the ideal \small
\[
I^* := \left\langle \begin{array}{c} x_3 x_5, \,x_0 x_5, \,x_2 x_3,
\,x_1 x_2-x_5^2, \,x_0 x_2, \,x_0
x_1-x_3^2, \,x_1^2 \beta \,+x_1 x_4 \beta \,+x_1 x_4+x_4^2, \, \\
x_1 x_5^2 \beta \,+x_4 x_5^2 \beta \,+x_2 x_4^2+x_4 x_5^2, \,x_1
x_3^2 \beta \,+x_3^2 x_4
\beta \,+x_3^2 x_4+x_0 x_4^2, \, \\
x_2 x_4 x_5^2 \beta \,+x_5^4 \beta \,+x_2^2 x_4^2+x_2 x_4 x_5^2,
\,x_3^4 \beta \,+x_0 x_3^2 x_4 \beta \,+x_0 x_3^2 x_4+x_0^2 x_4^2
\end{array} \right\rangle
\]
\normalsize which is reduced and has the following four associated
primes: \small
\[
\begin{array}{l}
H_1:= \langle x_5, \,x_2, \,x_1+x_4, \,x_3^2+x_0 x_4 \rangle, \,
B_1:= \lan x_5, \,x_2, \,x_1 \beta \,+x_4, \,x_0 x_1-x_3^2,
\,x_3^2 \beta \,+x_0 x_4\ran, \,
\\
H_2 := \langle x_3, \, x_0, \, x_1+x_4, \,x_2 x_4+x_5^2\ran, \,
B_2 := \lan x_3, \,x_0, \,x_1 \beta \,+x_4, \,x_1 x_2-x_5^2,
\,x_5^2 \beta \,+x_2 x_4\ran.
\end{array}
\]
\normalsize A quick observation reveals that
\begin{enumerate}\item  $H_1$ and $B_1$ intersect in a
tacnode at $[1,0,0,0,0,0]$;

\item  $H_2$ and $B_2$ intersect in a
tacnode at $[0,0,1,0,0,0]$;
\item  $H_1$ and $H_2$ intersect in a
node at $[0,1,0,0,-1,0]$;
\item  $B_1$ and $B_2$ intersect in a
node at $[0,1,0,0,-\beta, 0]$;
\item $H_i$ and $B_j$ do not
intersect if $i \ne j$.
\end{enumerate}
We conclude that  $I^*$ defines a curve isomorphic to the conjoined
snowmen curve $C^\ast$.

\end{proof}

\begin{prop}  $C^\flat$ is Chow (strictly) semistable.
\end{prop}
\begin{proof} We first show that
 $C$ is Chow strictly semistable with respect to $\rho_1$.
    Consider the normalization $\nu : E \to C$, $\nu\inv(p) =
\{q, r\}$.  We have $ord_q(\nu^*x_i) + r_i \ge 2$ and
$ord_r(\nu^*x_i) + r_i \ge 2$ for all $i$ where $r_i$'s denote the
weights of $\rho_1$, and it follows from Lemma 1.4 of \cite{Sch}
that
\[
e_{\rho_1}(C) = e_{\rho_1}(E)_q + e_{\rho_1}(E)_r \ge 2^2 + 2^2 =
8
\]
where $e_{\rho_1}(C)$ is the Hilbert-Samuel multiplicity of $C$
with respect to $\rho_1$.
 On the other hand, the right hand side of
the inequality in Theorem 1.1 of \cite{Sch} is $\frac23 \cdot 8
\cdot (2 + 1) = 8$. Hence $C$ is either Chow strictly semistable
or unstable with respect to $\rho_1$. It cannot be the latter
since $C$ is Hilbert semistable, and we conclude that $C$ (and
hence $C^\flat$) is Chow strictly semistable with respect to
$\rho_1$.

  By Proposition~\ref{P:link}, $C$ is Chow
semistable. Since $C$ is in the basin of attraction of $C^\flat$
(with respect to $\rho_1$), it follows that $C^\flat$ is also Chow
semistable and $C$ and $C^\flat$ are identified in the GIT
quotient space.
\end{proof}

\begin{prop} $C^\ast$ is Chow semistable.
\end{prop}

\begin{proof} We first show that
 $C^\ast$ is Chow strictly semistable with respect to $\rho_2$. We
 start with the ideal (\ref{E:Cflat}) of $C^\flat$.
 Let $\rho_2$ be the one-parameter subgroups with weights
$(0,2,2,1,2,2)$. We shall use Proposition~\ref{P:algorithm} to
show that $C^\flat$ is Hilbert unstable with respect to $\rho_2$:
The degree two monomials not in the initial ideal are
\[
x_0^2, \,x_0 x_3, \,x_0 x_4, \,x_1 x_3, \,x_1 x_4, \,x_1 x_5,
\,x_2^2, \,x_2 x_4, \,x_2 x_5, \,x_3^2, \,x_3 x_4, \,x_4^2, \,x_4
x_5, \,x_5^2
\]
whose $\rho_2$ weights sum up to 43. On the other hand, the
average weight is $\frac{2\cdot P(2)}{6} \cdot (2 + 2 + 1 + 2 + 2)
= 42$. Hence the Hilbert-Mumford index $\mu_2$ of the second
Hilbert point of $C^\flat$ with respect to $\rho_2$ is $6 \cdot
(42 - 43) = -6$. The degree three monomials not in the initial
ideal are
\[
\begin{array}{c}
x_0^3, \,x_0^2 x_3, \,x_0^2 x_4, \,x_0 x_3^2, \,x_0 x_3 x_4, \,x_1
x_3^2, \,x_1 x_3 x_4, \,x_1 x_4^2, \,x_1 x_4 x_5, \,x_2^3, \,
\\
x_2^2 x_4, \, x_2^2 x_5, \,x_2 x_4^2, \,x_2 x_4 x_5, \,x_2 x_5^2,
\,x_3^3, \,x_3^2 x_4, \,x_3 x_4^2, \,x_4^3, \,x_4^2 x_5, \,x_4
x_5^2, \,x_5^3
\end{array}
\]
whose $\rho_2$ weights sum up to 101, and the average weight is
$\frac{3\cdot P(3)}{6} \cdot (2 + 2 + 1 + 2 + 2) = 99$. Hence we
have $\mu_3 := \mu([C^\flat]_3, \rho_2) = 6 \cdot (99-101) = -12$.
Since $\mu_3 = 2 \mu_2 < 0$, it follows from
Proposition~\ref{P:all-m} that $C^\flat$ is $m$-Hilbert unstable
with respect to $\rho_2$ for all $m \ge 2$. By
Proposition~\ref{P:link} we conclude that $C^\flat$ is Chow
strictly semistable or unstable with respect to $\rho_2$. But the
latter cannot be the case since $C^\flat$ is Chow semistable.

  Since $C^\flat$ is in the basin of attraction of
$C^\ast$ (with respect to $\rho_2$), it follows that $C^\ast$ is
also Chow semistable and $C^\flat$ and $C^\ast$ are identified in
the GIT quotient space.

\end{proof}
Since
elliptic bridges and tacnodal curves belong to the basin of
attraction of the conjoined snowmen (Proposition~\ref{P:eb}), we
have

\begin{coro}\label{C:Zcs} All elliptic bridges
and tacnodal curves are identified in $\Mcs$.
\end{coro}
\noindent This completes the proof of Chow semistability of c-stable curves,
the second part of Theorem~\ref{T:main2}.
%\begin{figure}[h]
%\centerline{\scalebox{0.5}{\psfig{figure=cateye1.eps}}} \caption{A
%cat-eye}\label{F:cat-eye}
%\end{figure}

%\begin{figure}[h]
%\centerline{\scalebox{0.5}{\psfig{figure=ox.eps}}} \caption{An
%ox}\label{F:ox}
%\end{figure}

\section{Modular interpretation of the log canonical models}

\subsection{The first divisorial
contraction}\label{S:divisorial}

In this section, we summarize the result of \cite{HH1}.  In his
thesis \cite{Sch}, David Schubert considered the Chow stability of
tri-canonical curves and proved that a tri-canonical curve of
genus $\ge 3$ has a GIT stable Chow point if and only if it is
pseudostable. The corresponding statement was proved for genus two
curves by the authors in \cite{HL1}. A complete curve is {\em
pseudostable} if
\begin{enumerate}
\item it is connected, reduced, and
has only nodes and cusps as singularities;
 \item every subcurve of genus one meets the rest of
the curve in at least two points;
\item the canonical sheaf of the curve is ample.
\end{enumerate}
The last condition means that each subcurve of genus zero meets the
rest of the curve in at least three points.

Schubert also proved that there is no strictly semistable points
and constructed the moduli space $\Mps$ of pseudostable curves.
Our minimal model program is guided by the log canonical divisor
$K_{\bar{\cM}_3} + \a \d$ which is ample for $ 9/11 < \a \le 1$
but contracts a unique extremal ray consisting of elliptic tails
at $\a = 9/11$. In \cite{HH1}  is shown that the result of this
contraction, the log canonical model $\M_g(9/11)$, is isomorphic
to the moduli space $\M_g^{ps}$ for $g \ge 3$. For $g =2 $, this
was established by the authors in \cite{HL1}. More precisely,

\begin{thma} {\rm (1)}\, There is a natural morphism of stacks
 $\mathcal T : \FMg \to
\FMps$ that replaces each elliptic tail with a cusp, which
descends to a birational contraction $T : \Mg \to \M_g^{ps}$ with
exceptional locus $\D_1$;

\noindent {\rm (2)}\, The log canonical model $\Mg(9/11)$ is
isomorphic to $\M_g^{ps}$.
\end{thma}
 Let $\d^{ps}$ denote the divisor of singular curves in $\cMps$
and let $\la^{ps}$ denote the divisor class of $\bar \cM^{ps}_3$
that agrees with the Hodge class $\la$ on the open substack
$\cM_3$ of smooth curves. Since the locus of cuspidal curves is of
codimension two, $\la^{ps}$ is uniquely determined. Also, since
$T$ is a divisorial contraction, $K_{\cMps}$ is $\Q$-Cartier and
the equation $K_{\bar \cM_3} = 13 \la - 2 \d$ on the open substack
$\cM_3$ extends to
\[
K_{\cMps} = 13 \la^{ps} - 2 \d^{ps}
\]
on $\cMps$.

\begin{figure}[!t]
  \begin{center}
    \includegraphics[width=6.5in]{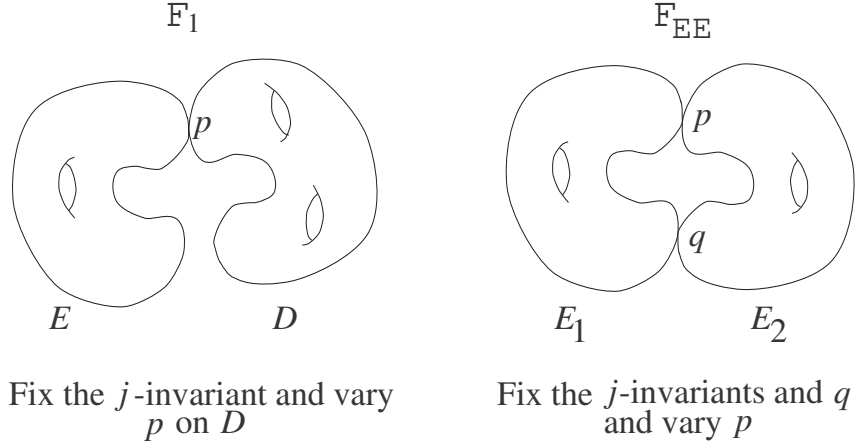}
  \end{center}

  \caption{$F_1$ and $F_{EE}$}
  \label{F:testcurves}
\end{figure}

\begin{lemma} $K^{ps}(\a) :=
K_{\cMps} + \a \d^{ps}$ is ample and $\M_3(\a) \simeq \M_3(9/11)$
for $7/10 < \a \le 9/11$.

\end{lemma}
\begin{proof} Since $\M_3^{ps} \simeq \M_3(9/11)$, $K^{ps}(9/11)$
is ample on $\M_3^{ps}$. Also,  $K^{ps}(7/10)$ is nef on
$\M_3^{ps}$ due to Faber (\cite{F1}): The cone $\nef^1(\M_3)$ of
nef divisors is generated by $\la$, $12 \la - \d_0$ and $10\la -
\d_0 - 2\d_1$, and $T^*( K_{\cMps} + 7/10 \, \d^{ps}) =
\frac{13}{10} (10 \la - \d_0 - 2\d_1)$.

 We can write $K^{ps}(\a)$ as a linear combination
 of $K^{ps}(9/11)$ and $K^{ps}(7/10)$:
\[
K^{ps}(\a) = \frac{11}{13}(10 \a - 7) K^{ps}(9/11) + \frac{10}{13}(9
- 11\a) K^{ps}(7/10).
\]
Since the first coefficient is positive and the second is
nonnegative for $7/10 < \a \le 9/11$, the ampleness follows.
\end{proof}
From Faber's result, one also readily deduces that the nef cone of
$\M_3^{ps}$ is generated by $10\la^{ps} - \d^{ps}$ and $12\la^{ps}
- \d^{ps}$: Since $\M_3^{ps}$ is of Picard number two, it suffices
to show that these are nef and extremal. But these divisors pull
back to the nef generators $10\la - \d_0 - 2\d_1$ and $12\la -
\d_0$, and contract the loci $T(F_{EE})$ and $T(F_1)$
(Figure~\ref{F:testcurves}), respectively.
%\begin{figure}[htb]
%\centerline{\scalebox{0.6}{\psfig{figure=testcurves.eps}}}
%\caption{$F_1$ and $F_{EE}$}\label{F:testcurves}
%\end{figure}

\subsection{The small
contraction}\label{S:lcmodel-Chow} In this section, we prove the
part (1) of Theorem~\ref{T:main1}:

\

\ndt {\bf Theorem 1.} {\rm (1)}\, {\it
 $\M_3(7/10)$ is isomorphic to the moduli space
$\M_3^{cs}$ of c-stable curves of genus three. }

 \

The moduli space  $\M_3^{cs}$ of c-stable curves (resp.
$\M_3^{hs}$ of h-stable curves)  was constructed in \S\ref{S:GIT}
as a GIT quotient of the Chow variety $\Ch$ (resp. the Hilbert
scheme $\H$) of the bicanonical curves of genus three. We shall
employ the following strategy to prove Theorem~\ref{T:main1}.(1):
We first show the existence of the log canonical model
$\M_3(7/10)$ by applying Kawamata-Shokurov base point freeness
theorem on $K_{\cMps} + 7/10 \d^{ps}$. Then we'll show that the
polarizations on $\M_3(7/10)$ and $\M_3^{cs}$ agree when pulled
back to $\M_3^{ps}$. Since the two varieties are normal and
isomorphic away from loci of codimension $\ge 2$, the theorem
follows by Hartog's Lemma.

\

\noindent Recall the base point freeness theorem:

\begin{thma}(see, e.g. \cite{KMM}, \cite{KM})
  Let $(X, \Delta)$ be a proper klt pair with
$\D$ effective. Let $D$ be a Cartier divisor such that $a D - K_X -
\D$ is nef and big for some $a > 0$. Then $|b D|$ has no basepoint
for all $b \gg 0$.
\end{thma}
Let us apply this theorem to the pair $(X, \Delta) = (\M_3^{ps},
\a\dps)$ for some $\a < 7/10$, $D = K_{\cMps} + 7/10 \, \dps$ and
$a = 2$.

\begin{prop}\label{P:morphism}
 The linear system $\left|m(K_{\cMps} + 7/10 \d^{ps})\right|$
is base point free for sufficiently large and divisible $m$, and
 the associated morphism $\Psi : \M_3^{ps} \to
\M_3(7/10)$  is a small contraction.
\end{prop}

\begin{proof}
Due to the base point freeness theorem, it suffices to establish
that
\begin{enumerate}
\item[(A)] $(\M_3^{ps}, \a \dps)$ is klt for $\a < 1$;

\item[(B)] $2 (K_{\cMps} + 7/10 \, \dps) - (K_{\cMps} + \a \dps)$ is nef and big
for $\a < 7/10$.

\end{enumerate}

\ndt (B) follows immediately from that $K_{\cMps} + \b\dps$ is ample
for $7/10 < \b \le 9/11$. For (A), recall the log discrepancy
formula
\begin{equation}\label{E:logdisc}
K_{\bar \cM_3} + \a \d = T^*(K_{\cMps} + \a \dps) +
  (9 - 11\a) \d_1.
\end{equation}
Since $\bar \cM_3$ is smooth, (\ref{E:logdisc}) implies that the
stacky pair $(\cMps, \a\dps)$ is klt for $9 - 11\a \ge -1$ and $\a
< 1$. That $(\cMps, \a\dps)$ is klt now follows from \cite{KM},
5.20. (See also \cite{HH1}, A.13.)

We now prove that $\Psi$ is a birational morphism. Let $C \subset
\M_3$ be a curve that is not contained in the boundary. The
Moriwaki divisor
\[
A := 28 \la - 3\d_0 - 8\d_1
\]
intersects with any such curve non-negatively \cite{Mor}. But
$K_{\cMps} + 7/10 \d^{ps}$ pulls back to a line bundle
proportional to $10\la -\d_0-2\d_1$, which in turn is a positive
rational multiple of $A + 2 \la + 2 \d_1$ and
\begin{equation}\label{E:Moriwaki}
(A + 2 \la + 2 \d_1).C \ge 2 \la.C > 0
\end{equation}
where the last inequality follows since $\la$ gives rise to the
Torelli map.

It remains to show that $\Psi$ is a small contraction. If $\Psi$
is a divisorial contraction, it contracts an irreducible divisor.
Since  (\ref{E:Moriwaki}) implies that $\Psi$ does not contract
the hyperelliptic divisor, it must contract $\d^{ps}$. This would
force $T^*(K_{\cMps} + 7/10 \d^{ps})$ to contract both $\d_0$ and
$\d_1$, but a divisor contracting $\d_0$ and $\d_1$ must be
proportional to $\lambda$.

\end{proof}

\begin{proof}[Proof of Theorem~\ref{T:main1}, part (1)]
Let $Z^{cs} \subset \M_3^{cs}$ (resp. $Z^+ \subset \M_3^{hs}$)
denote the locus  consisting of tacnodal c-stable curves and
elliptic bridges (resp. h-stable curves with tacnodes). $Z^{cs}$
is a point by Corollary~\ref{C:Zcs}. By definition, a c-stable
curve without a tacnode or an elliptic bridge is pseudo-stable;
Similarly, an h-stable curve without a tacnode is pseudo-stable.
Hence there are isomorphisms
\begin{equation}\label{E:isoms}
\M_3^{ps}\setminus Z \simeq \M_3^{cs}\setminus Z^{cs} \simeq
\M_3^{hs}\setminus Z^+
\end{equation}
Note that the loci $Z$, $Z^{cs}$ and $Z^+$ are of codimension $\ge
2$.

Let $\pi : \H \to \Mhs$ denote the quotient map. The projective
structure on $\Mhs$ is given by the invariant sections of
$\SS_{\H}(+m)$ on $\H$. Hence $\pi^*(\O_{\Mhs}(1)) = \O_{\H}(1)$ which
is proportional to (\cite{Vie}, \cite{HH2})
\begin{equation}\label{E:poln}
\left( 10 - \frac 3{2m}\right) \la - \d.
\end{equation}
 Hence for large enough $m$,  $\left(
10 - \frac 3{2m}\right) \la - \d$ descends to an ample
$\Q$-divisor on $\Mhs$. It follows that $\la$ (and consequently,
$\d$)   descends to a Cartier divisor since it is
 proportional to the
difference of $\left(10-\frac3{2m}\right)\la - \d$ and
$\left(10-\frac3{2(m+1)}\right)\la - \d$. Let $\lhs$ and $\dhs$
denote the Cartier divisors on $\Mhs$ that pull back to $\la$ and
$\d$.

The cycle map $\cycle : \H \to \Ch$  is $\SL(6)$-equivariant and
descends to
\begin{equation}\label{E:Psiplus}
\Psi^+ : \M_3^{hs} \to \M_3^{cs}.
\end{equation}
By a theorem of Mumford \cite{Mum}, the polarization on
$\M_3^{cs}$ pulls back by $\Psi^+$ to $10\lhs - \dhs$ on
$\M_3^{hs}$. This line bundle is isomorphic to a positive rational
multiple of $K_{\cMps} + (7/10)\dps$ over $\M_3^{ps} \setminus Z$,
via the isomorphism (\ref{E:isoms}). The varieties are normal, and
a section of a line bundle over a normal variety defined except
over a codimension $\ge 2$ locus extends uniquely to a global
section. We conclude that
\[
\G(\M_3^{ps}, m (K_{\cMps} + (7/10)\dps)) \simeq \G(\M_3^{cs}, \mcl
O_{\M_3^{cs}}(+m))
\]
which leads to
\[
\begin{array}{ll}
\M_3(7/10) & \simeq \Proj \ds_{m\in\Z_+}\G(\M_3^{ps}, m (K_{\cMps} +
(7/10)\dps)) \\
& \simeq  \Proj \ds_{m\in\Z_+} \G(\M_3^{cs}, \mcl O_{\M_3^{cs}}(+m))
\simeq \M_3^{cs}
\end{array}
\]
where the first isomorphism follows from the projection formula.

It remains to prove that the exceptional locus of $\Psi$ is $Z$.
Recall the test curve $F_{EE}$  from
Section~\ref{S:divisorial} (Figure~\ref{F:testcurves}). Since there is a $T(F_{EE})$-curve
that passes through a general point of $Z$ and $10\la-\d_0-2\d_1 =
T^*(10\la^{ps}-\d^{ps})$ contracts $F_{EE}$, it follows that the
exceptional loci of $\Psi$ contains $Z$. We showed  that
$K_{\cMps} + (7/10)\dps$ defines a small contraction from $\Mps$
to $\Mcs$ (Proposition~\ref{P:morphism}), and observed above that
under the isomorphism $\Mcs\setminus Z^{cs} \simeq \Mps\setminus
Z$, the polarization on $\Mcs$ is proportional to $K_{\cMps} +
(7/10)\dps$ on $\Mps\setminus Z$. It follows that $\Psi$ is an
isomorphism on the open subset $\M_3^{ps} \setminus Z$. Hence
$Excep(\Psi) = Z$, and this completes the proof of the first part
of Theorem~1.
\end{proof}

 Now that we have modular interpretations of $\M_3(9/11)$ and
$\M_3(7/10)$, we can describe the small contraction $\Psi$ induced
by $K_{\cMps} + 7/10 \, \dps$ in a concrete manner. It is an
isomorphism over the $\Mps\setminus Z$.  What does $\Psi$ do to
curves in $Z$?  That is, if $C$ is an elliptic bridge, what is the
c-stable curve that corresponds to $\Psi(C)$? Since $\M_3(7/10)$
is the moduli space of c-stable curves and $C$ is c-stable, the
semistable replacement theorem implies that $\Psi(C)$ is $[C]$
itself. Moreover, we proved in \S\ref{S:GIT} that elliptic bridges
and tacnodal c-stable curves are identified with the two snowmen
curve $C^\ast$ (Proposition~\ref{P:eb-unstable}) in the GIT
quotient. We conclude that

\begin{prop} $\Psi$ collapses $Z$ to the point in $\M_3^{cs}$ that
represents tacnodal curves.
\end{prop}

\subsection{The flip}\label{S:lcmodel-Hilb} In this section, we shall prove

\medskip

\ndt {\bf Theorem 1.} {\rm (2)} \, \textit{
 For $\a \in (17/28, 7/10)$, the log canonical model $\M_3(\a)$
  is isomorphic to  $\Mhs \simeq \Hbar$.
  Moreover, the following diagram is a flip in the sense of Mori
theory:}
\begin{equation}\label{D:flip}
\xymatrix{  \Mps  \ar[dr]_-{\Psi} & &
\Mhs \ar[dl]^-{\Psi^+} \\
& \M_3^{cs} & \\
}
\end{equation}

\

\begin{proof} We first prove the assertion for
$\a \in (7/10 - \e, 7/10)$ for small enough $\e$. It will be
extended to all $\a\in(17/28, 7/10)$ in Corollaries
\ref{C:17/28-1} and \ref{C:17/28-2}.

The ampleness of (\ref{E:poln}) implies that
\begin{equation}\label{E:hg} \Mhs \simeq \Proj \oplus_{s \in \Z_+}
\Gamma\left(\Mhs, s \left( \left(10-\frac3{2m}\right)\la -
\d\right)\right), \quad m \gg 0.
\end{equation}
Recall from \S\ref{S:lcmodel-Chow} that  $(\Mhs)^o :=
\Mhs\setminus Z^+$ is isomorphic to $\M_3^{ps}\setminus Z$. On
$\Mps$, $(10 - \frac 3{2m})\la - \d$ is proportional to $K_{\Mps}
+ \frac{14m-6}{20m-3} \, \dps$.  Since $\Mhs$ and $\M_3^{ps}$ are
normal and isomorphic away from the codimension $\ge 2$ loci $Z$
and $Z^+$, we have: \smaller
\begin{equation}\label{E:h-ps}
\bigoplus_{s \in \Z_+} \Gamma\left(\Mhs, s \left(
\left(10-\frac3{2m}\right)\la - \d\right)\right) = \bigoplus_{s
\in \Z_+} \Gamma\left(\Mps, s \left( \left( K_{\Mps} + (7/10 -
\e(m))\dps \right)\right)\right)
\end{equation}
\normalsize where $ \e(m) := 39/(200m - 30).$
From the log discrepancy formula
$
K_{\bar \cM_3} + \a \d = T^*(K_{\Mps} + \a \dps) + (9-11\a)\d_1$ for
$T : \M_3 \to \Mps$,
we get
\begin{equation}\label{E:ps-mg}
\begin{array}{rl}
 \G(\M_3, s(K_{\bar \cM_3} + \a\d) ) & = \G(\M_3, s(T^*(K_{\Mps} + \a \dps)
+ (9-11\a)\dps))\\
& = \G(\Mps, s( K_{\Mps} + \a \dps))
\end{array}
\end{equation}
for all $\a$. Combining (\ref{E:hg})$\sim$(\ref{E:ps-mg}) gives the
desired isomorphism
\[
\begin{array}{rl}
\Mhs & \simeq \Proj \ds_{s\in \Z} \G(\M_3, s(K_{\bar \cM_3} +
(7/10-\e(m))\d) ) \\
& = \M_3(7/10 - \e(m)).
\end{array}
\]
Recall that $\Psi$ is the birational contraction with exceptional
locus $Z$ which is the extremal face for the divisor $K_{\Mps} +
(7/10 - \e) \dps$ for a small enough $\e$. Hence to show that
(\ref{D:flip}) is a Mori flip, we need to establish that

\begin{enumerate}
\item $\Psi^+$ is a small contraction;
\item The strict transformation of $K_{\Mps}+ (7/10 - \e)\dps$ is
$\Q$-Cartier and $\Psi^+$-ample.
\end{enumerate}
The first item follows since $\Psi^+$ has exceptional locus $Z^+$
is of dimension one (Corollary~\ref{C:Z^+}). The strict
transformation of $K_{\Mps}+ (7/10 - \e)\dps$ is $K_{\Mhs} + (7/10
- \e)\dhs$. Note that $K_{\Mhs}$ agrees with $13\lhs - 2 \dhs$ on
the open set $\Mhs\setminus Z^+$ since $K_{\Mps} = 13\la^{ps} -
2\dps$. This implies that $K_{\Mhs} = 13\lhs - 2 \dhs$ on $\Mhs$,
and $K_{\Mhs} + (7/10 - \e)\dhs = 13\lhs - (13/10 + \e) \dhs$ is
$\Q$-Cartier. It is also $\Psi^+$-ample, being a positive rational
multiple of the polarization (\ref{E:poln}).

\end{proof}

\begin{prop}\label{P:17/28}
 $K_{\M_3^{hs}}+\frac{17}{28}\dhs$ is nef and big on $\M_3^{hs}$ and
 has a unique extremal ray generated by hyperelliptic curves.
\end{prop}

\begin{proof}  Since $K_{\Mhs} = 13 \lhs - 2\dhs$,
$K_{\M_3^{hs}}+\frac{17}{28}\dhs$ is proportional to
$28\lhs-3\dhs$. On the other hand, the Moriwaki divisor $A = 28\la
- 3\d_0 - 8\d_1$ properly transforms to $28\lhs - 3\dhs$ on
$\Mhs$. Therefore, $K_{\M_3^{hs}}+\frac{17}{28}\dhs$ is big. Since
the hyperelliptic locus is equal to $9 \la - \d_0 - 3 \d_1$ on
$\M_3$, $h = 9 \lhs - \dhs$ on $\Mhs$ and we have
\[
28\lhs-3\dhs=(10\lhs-\dhs)+2(9\lhs-\dhs) = (10\lhs-\dhs) + 2h.
\]
The divisor $10\lhs-\dhs$ is nef, since it is proportional to
$K_{\M_3^{hs}}+\frac{7}{10}\dhs$ which is a limit of ample
divisors. It follows that to show the nefness of $28\lhs-3\dhs$,
it suffices to prove that the divisor non-negatively intersects
with curves in $h$. In fact, we claim that $28\lhs-3\dhs$
is trivial on $h$.

Let $B_i \subset \M_3$ denote the locus of curves obtained by taking
the stabilization of the admissible cover of $C \in \M_{0,8}$
consisting of  smooth rational curves $C_1$ with $i$ marked points
and $C_2$ with $8-i$ marked points meeting in one node. Abusing
notation, let $h$ denote the locus in $\M_3$ of hyperelliptic
curves. For $g = 3$, we have three boundary divisors of $h$ which
generate the rational Picard group:
\begin{enumerate}
\item $B_2$ consists of irreducible curves with one node;
\item $B_3$ consists of elliptic tails;
\item $B_4$ consists of elliptic bridges.
\end{enumerate}
Among these, $B_3$ is contracted by $T$ and $B_4$, by $\Psi$. Hence
the rational Picard group of the hyperelliptic locus $h^{cs}$ of
$\Mcs$ is generated by the image of $B_2$. The small contraction
$\Psi^+ : \Mhs \to \Mcs$ restricts to a small contraction  on $h$,
and $\Psi^+|_h$ induces an isomorphism $\Pic(h)\ten \Q \simeq
\Pic(h^{cs})\ten \Q$ of the Picard groups. We conclude that
$\Pic(h)\ten \Q$ is generated by the image of $B_2$.

We summarize some results from \S2.1 of \cite{R}: Any smooth
hyperelliptic curve of genus three is a divisor of type $(2,4)$ in
$\P^1\times \P^1$. Let $F_h$ denote a pencil of these divisors,
which is equivalent to twice the curve class $F_h'$ in $B_2$
obtained by letting one of the six marked points move. Since $F_h'$
is in $\bar{\mathcal M}_3\setminus (\d_1\cup \{ \textrm{elliptic
bridges}\})$, the following intersection computation on
$\bar{\mathcal M}_3$ carries over to $\Mhs$:
\[
F_h'.\d_0 = 14, F_h'.\d_1 = 0, F_h'.\la = 3/2.
\]
It follows that $(28\lhs-3\dhs)\cdot F_h'=0$. Since
$\Pic(h)\otimes\Q$ is generated by the image of $B_2$, the divisor
$28\lhs-3\dhs$ is trivial on $h$. Therefore, $28\lhs-3\dhs$ is nef
and it contracts $h$ to a point.

Since the Mori cone of $\M_3^{hs}$ is of dimension two and
$K_{\M_3^{hs}}+\frac{17}{28}\dhs$ is big, the unicity of the
extremal ray is obvious.

\end{proof}

\begin{coro}\label{C:17/28-1}
$K_{\M_3^{hs}}+\a\dhs$ is ample if $\a\in (\frac{17}{28}, \
\frac{7}{10})\cap\Q$.
\end{coro}

\begin{proof} Given $\a\in
(\frac{17}{28}, \ \frac{7}{10})\cap\Q$ and small $\e$,
$K_{\M_3^{hs}}+\a\dhs$ is a positive multiple of the linear
combination
\[
\left(\a-17/28\right)\left(K_{\M_3^{hs}}
+(7/10-\epsilon)\dhs\right) + (7/10 -\a -
\e)\left(K_{\M_3^{hs}}+(17/28)\dhs\right).
\]
%\[\frac1{{\frac{7}{10} - \frac{17}{28}
%- \e}}\left[\left(\a-17/28\right)\left(K_{\M_3^{hs}}
%+(\frac{7}{10}-\epsilon)\dhs\right) + (7/10 -\a -
%\e)\left(K_{\M_3^{hs}}+\frac{17}{28}\dhs\right)\right].\]
Since the divisor $K_{\M_3^{hs}}+(7/10-\epsilon)\dhs$ is ample for
small enough $\e$ and $K_{\M_3^{hs}}+(17/28)\dhs$ is nef,
 $K_{\M_3^{hs}}+\a\dhs$ is ample.
\end{proof}

We have established that
\begin{coro}\label{C:17/28-2} For $\a \in (17/28, 7/10)$, $\M_3(\a)$ is
isomorphic to $\Mhs$.
\end{coro}
This completes the proof of the second part of
Theorem~\ref{T:main1}.

\subsection{The second divisorial
contraction}\label{S:17/28}

This is the last (nontrivial) step in the Mori program for $\M_3$.
We have shown in Proposition~\ref{P:17/28} that $K_{\Mhs} + (17/28) \dhs$
contracts the hyperelliptic locus, and we aim to describe the
resulting log canonical model $\M_3(17/28)$.

We first need to show that $\M_3(17/28)$ exists, and we use the
base point freeness theorem as we did in
Proposition~\ref{P:morphism}. In general, if $(X, \D)$ is klt then
so is $(X, a \D)$ for any $0 < a < 1$. Since $(\Mhs, \a \dhs)$ is
klt for any $\a \in (17/28, 7/10)$, $(\Mhs, a \dhs)$ is klt for
any $0 < a < 7/10$.

Choose $\e$ such that $7/10 - 17/28 < \e < 2 (7/10 - 17/28)$. We
have
\begin{equation}\label{E:bpf-mhs}
2 ( K_{\Mhs} + (17/28) \dhs ) - (K_{\Mhs} + (7/10 - \e)\dhs) =
K_{\Mhs} + \beta \dhs
\end{equation}
where $\beta := 2\cdot 17/28 - 7/10 + \e$. The log canonical
divisor (\ref{E:bpf-mhs}) is nef and big since $17/28 < \beta <
7/10$. In fact, it is ample. The base point freeness theorem
implies that $|b( K_{\Mhs} + (17/28)\dhs)|$ is base point free for
$b \gg 0$. Hence
\begin{lemma} $\M_3(17/28)$ exists as a projective variety, and
there is a divisorial contraction
\[
\Theta : \Mhs \to \M_3(17/28)
\]
with exceptional locus $h$.
\end{lemma}

%Note that  the canonical divisor is not very ample on $C$ if $C$
%has an elliptic tail or an elliptic bridge or a tacnode: One can
%verify the first two just by restricting $\o_C$ to an elliptic
%tail or bridge. To see the tacnode case, let $p$ denote the
%tacnode and consider $H^1(C, \mfk m_p^2 \o_C)$. We have the
%isomorphisms
%$$ H^1(C, \mfk m_p^2 \o_C)^* \simeq Hom(\mfk m_p^2 \o_C, \o_C) \simeq
%H^0(\til{C}, \mcl O_{\tilde{C}}(p_1+p_2)) \ne 0$$ where $\tilde{C}
%\to C$ is the partial normalization resolving $p$, and $p_1$,
%$p_2$ are the points sitting over $p$.

We recall some classical results on the GIT of plane quartics \cite{GIT},
\cite{Art}.
Consider the natural action of $\PGL(6)$ on the space
$\P(\Gamma(\cO_{\P^2}(4)))$ of plane quartics. With respect to
this action, a plane quartic curve $C$ is
\begin{enumerate}
%\item is unstable if and only if it is a cubic with an inflectional tangent
%line or has triple points;
%
\item stable if and only if it has ordinary nodes and
cusps as singularity;

\item strictly semistable if it is a double conic or has a tacnode.
Moreover, $C$ belongs to a minimal orbit if and only if it is
either a double conic or the union of two tangent conics (where at
least one is smooth).
\end{enumerate}
The minimal orbit statement in (2) implies that in the GIT
quotient space $\bar Q$, an irreducible tacnodal curve is
identified with the corresponding cat-eye, as in $\Mhs$.

\

\ndt {\bf Theorem 1.} {\rm (3)} \,  \textit{$\M_3(17/28)$ is
isomorphic to the compact space of plane quartics $\bar Q:=
\P(\Gamma(\SS_{\P_2}(+4)))\mod \PGL(6)$.}

\smallskip

\begin{proof} Consider the universal quartic curve $\mcl X$ over
$\P := \P(\Gamma(\cO_{\P^2}(4)))^{ss}$:
\[
\xymatrix{ \mcl X \ar@{^(->}[r] \ar[d] & \P \ex \P^2 \ar@{^(->}[r]
& \P \ex v_2(\P^2) \subset \P \ex \P^5 \\
\P & & \\
}
\]
where $v_2$ denotes the second Veronese embedding. By abusing
notation, let $\mcl X$ denote the image of $\mcl X$ in $\P \ex
\P^5$.  Let $y\in \bar Q$ denote the point corresponding to the
double conic. Away from the orbit of the double conic $\mcl X_y$,
$\mcl X \inj \P\ex\P^5$ is a family of bicanonical h-semistable
curves, and induces a map from $\P\setminus O(\mcl X_y) $ to $
\H^{ss}$ and subsequently to the quotient $\H\mod\SL$. This map is
$\PGL(3)$ invariant, and it descends to give $f : \bar
Q\setminus\{y\} \to \Mhs$. Let $Z^q \subset \bar Q$ denote the
locus of tacnodal curves. Over the stable locus, we have
isomorphisms
$$\bar Q\setminus (\{y\}\cup Z^q) \simeq \Mhs\setminus (Z^+\cup
h) \simeq \M_3\setminus (\d_1\cup \{\textrm{elliptic bridges}\}
\cup h)$$ where $ f$ induces the first isomorphism. The loci $Z^q$
and $Z^+$ are of codimension $\ge 2$, and $f|_{Z^q}$ is bijective
(see item (2) above and the subsequent remark). It follows that
the inverse rational map $f^{-1}$ is regular on $Z^+$, giving an
isomorphism $$\bar Q\setminus \{y\} \simeq \Mhs\setminus h \simeq
\M_3(17/28) \setminus \{\Theta(h)\}.$$ The assertion now follows
by applying Hartog's Lemma again.
\end{proof}
Let $\la^{\bar Q}$ and $\,  \d^{\bar Q}$ denote the unique divisor
extending $\lhs$ and $\dhs$ on the open set $\bar Q\setminus \{y\}
\simeq \Mhs\setminus Z^+$.
\begin{lemma} $K_{\bar Q} + \a \,  \d^{\bar Q}$ is ample on $\bar Q$ for
$\a \in (5/9, 17/28]$ and $K_{\bar Q} + (5/9) \,  \d^{\bar Q}$ is
trivial.
\end{lemma}

\begin{proof}  The identities $K_{\Mhs} = 13 \lhs
- 2\dhs$ and $h = 9 \lhs - \dhs$ carry over to $\bar Q$ and give
\[
K_{\bar Q} + 5/9 \,  \d^{\bar Q} = 13\la^{\bar Q} - 2\,  \d^{\bar
Q} + 5/9 \,  \d^{\bar Q} = (13/9 - 2 + 5/9) \,  \d^{\bar Q} = 0.
\]
Hence $K_{\bar Q} + 5/9 \,  \d^{\bar Q}$ is trivial. From this
follows that $K_{\bar Q} + \a \,  \d^{\bar Q}$ is ample for $\a \in
(5/9, 17/28]$ since it is a linear combination
\[
a \left(K_{\bar Q} + 5/9 \,  \d^{\bar Q}\right) +  b \left(K_{\bar
Q} + 17/28 \,  \d^{\bar Q}\right)
\]
for some positive rational numbers $a, b$ (determined by $\a$).
\end{proof}
\begin{coro} $\M_3(5/9)$ is a point.
\end{coro}

\section{Relation to other moduli spaces: Work of Hassett, Hacking and
Kondo}

There are various constructions of compact moduli spaces of plane
quartics \cite{Has1}, \cite{Hac}, \cite{Kon}. How do these moduli
spaces fit in our minimal model program? In this section we give a
brief sketch of these moduli spaces and show that they are indeed
log canonical models for $\M_3$.

\subsection{Hassett's moduli space $\HP$}\label{S:Has}
In \cite{Has1}, B. Hassett constructed a compact moduli space
$\bar{\mathcal P}_4$  by taking the closure of the $\mcl P_4$ in
the connected moduli scheme of {\it smoothable} stable log
surfaces, where $\mcl P_4$ is the quasi-projective GIT moduli
space of smooth plane quartics. A stable log surface is a  pair
consisting of a surface $S$ and a curve $C\subset S$ such that
$(S, C)$ has semi-log canonical singularities and $K_S+C$ is
ample. It is said to be smoothable if there is a a one parameter
family of deformations $(\mcl S, \mcl C)$ whose general fiber is a
smooth pair consisting of $\P^2$ with a smooth plane quartic curve and
both $\Q$-divisors $K_{\mcl S}+\mcl C$ and $\mcl C$ are
$\Q$-Cartier.

There is a forgetting morphism $F: \HP \to \M_3$ defined by $F((S,
C))= C$ which in fact, is isomorphism: For each curve $C$ in $\M_3$,
Hassett explicitly constructs the unique corresponding surface $S$
with $(S, C)\in \HP$. Then he shows that the morphism $F$ is proper,
birational, and locally an isomorphism.

The cone of effective divisors $\bar{NE}^1(\bar{\mcl M}_3)$ is
generated by $\d_0, \d_1$ and $h$. A general element in each of
these divisors corresponds to the following stable log surface in
$\HP$:
\begin{enumerate}
\item $C$ in $\d_0$ $\Longleftrightarrow$ $(\P^2, C)$;
\item $C=C_1\cup_p C_2$ in $\d_1$ where $C_1$ and $C_2$ are
irreducible curves of genus two and one respectively
$\Longleftrightarrow$ $(S_1\cup_B S_2, C)$ where $S_1$ is the
toroidal blowup of $\P^2$ and $S_2=\P(1,2,3)$. The curve $C_1$ in
$S_1$ does not pass through the two singular points of type
$\frac{1}{2}(1,1)$ and $\frac{1}{3}(1,1)$ of the surface $S_1$. Also
the curve $C_2$ in $S_2$ does not pass through the two singular
points $\frac{1}{2}(1,1)$ and $\frac{1}{3}(1,2)$ of the surface
$S_2$. The curve $C=C_1\cup_p C_2$ meets the double curve $B$ at the
point $p$.
\item $C$ in $h$ $\Longleftrightarrow$ $(S, C)$ where
$S=\P(1,1,4)$. Since $C$ is a smooth hyperelliptic curve of genus
three, $C$ is regarded as a bisection of the rational surface
$\mathbb F_4=\P(\mcl O_{\P^1}\oplus\mcl O_{\P^1}(4))$. $S$ is
obtained by the contraction of the zero section. The curve $C$ does
not pass through the singular point of type $\frac{1}{4}(1,1)$.
\end{enumerate}
%Note that for any pairs $(S, C)\in\HP$, $C$ is a Cartier divisor and
%$6K_S$ is a Cartier divisor in $S$.

\subsection{Hacking's moduli space $\HP'$}\label{S:Hac}
Hacking gave an alternate  compactification of the moduli space of
plane curves of degree $d$ \cite{Hac} by employing a method
similar to \cite{Has1} but allowing worse singularities. We denote
his moduli space for $d = 4$ by $\HP'$ (His original notation
$\mcl M_3$ is unfortunately reserved for the moduli stack of
smooth curves of genus three.) A  {\it stable pair} of degree 4 is
a pair consisting of a surface $S$ and a curve $C\subset S$ such
that $(S, \beta C)$ has semi-log canonical singularities and
$K_S+\beta C$ is ample, where $\beta$ is a rational number
$\frac{3}{4}+\epsilon$ for a sufficient small positive number
$\e$. We also impose the condition that there be a one parameter
family of deformations $(\mcl S, \mcl C)$ whose general fiber is a
smooth pair $(\P^2, \textrm{smooth quartic curve})$. Both
$\Q$-divisors $K_{\mcl S} + \beta \mcl C$ and $\mcl C$ are
$\Q$-Cartier.

From his classification of stable surfaces of degree 4,
$\HP'=Z_0\cup Z_1\cup Z_2$ where $Z_1$ has codimension 1 and $Z_2$
has codimension 2 such that:
\begin{enumerate}
\item Any element in $Z_0$ is a pair $(\P^2, C)$ where $C$ is a
pseudo-stable plane curve of degree 4;
\item Any element in $Z_1$ is a pair $(S, C)$ where
$S=\P(1,1,4)$ and $C$ is a (degenerating) hyperelliptic curve of
genus 3. The curve $C$ does not pass through the singular point of
type $\frac{1}{4}(1,1)$;
\item Any element in $Z_2$ is a pair $(S, C)$ where $S=S_1\cup_B S_2$
is the union of two $\P(1,1,2)$s and $C$ is a (degenerating)
elliptic bridge. Both irreducible components $S_1$ and $S_2$ have
cyclic quotient singularities of type $\frac{1}{2}(1,1)$ on the
double curve $B=\P^1$. The curve $C$ does not pass through the
singular points of $S_1$ and $S_2$.
\end{enumerate}
If $(S, C)$ is a stable pair in $\HP'$ then $C$ has nodes and
cusps as singularities, and there is a forgetting morphism $F':
\HP'\to \M_3^{ps}$ defined by $F'((S, C))= C$.

\subsection{$\HP$ and $\HP'$ as log canonical models}
Hassett proved that $\HP$ is isomorphic to the moduli space
$\bar{M}_3$ of stable curves. Implicit in Hacking's work is that

\begin{prop} $\HP'$ is isomorphic to $\M_3(9/11)$.
\end{prop}
\begin{proof}
It is remarked in Hassett's paper that for $\beta > 5/6$, $K_S +
\beta C$ is ample for any pair $(S, C) \in \HP$. Therefore
$\HP(\beta) \simeq \HP$ for $\beta > 5/6$, where $\HP(\beta)$
denotes the moduli space of pairs constructed in \S\ref{S:Hac}
using the prescribed value of $\beta$ whereas $3/4 +\e$ was used
for $\beta$ in \S\ref{S:Hac}.

We consider what happens at $\beta = 5/6$. First, there is a
morphism $T' : \HP \to \HP(5/6)$ that associates  to $(S, C)$ the
pair $(S', C')$ constructed as follows: By assumption $(S, C)$ is
smoothable and there is a one-parameter family $(\mcl S, \mcl C)$,
$\pi : \mcl S \to Spec(k[[t]])$ whose
special fibre is $(S,C)$ and whose generic fibre is smooth. $S'$
is then the special fibre of the relative log canonical model
$\mcl S' := \Proj\ds_{m\ge 0} \pi_*(m (K_\mcl S + 5/6 \, \mcl C))$
and $\mcl C'$, the scheme theoretic image $f(\mcl C')$ where $f$
is the canonical fibration from $\mcl S$ to $\mcl S'$. Regardless
of the choice of the smoothing, $(S',C')$ is always the log
canonical model of $(S, C)$ (If $S$ has more than one component,
then $S'$ is the union of the log canonical models of the
components with the restriction of the log canonical divisor as
boundary divisor.)

 Since the
divisor $K_S + 5/6 \, C$ is ample for all pairs $(S,C) \in \HP$
except the ones with $S = (\textrm{toroidal blowup of $\P^2$}\cup
\P(1,2,3))$ and $C$ an elliptic tail, $T'$ is an isomorphism away
from the locus $\mcl D_1 \subset \HP$ of elliptic tails. For
elliptic tails $(S_1\cup_B S_2, C)$, restricting $K_S + 5/6 \, C$
on $S_2$, we find that
\[
(K_S+\frac{5}{6}C)|_{S_2}=K_{S_2}+\frac{5}{6}C+B=\cO_{S_2}(-6+5+1)=\cO_{S_2}.
\]
This means that $f : S \to S'$ contracts the elliptic component:
In fact, Hassett proves that $f : C \to C'$ replaces the elliptic
tail by an ordinary cusp. Also, the log canonical models of $(S,C)
\in \HP$ with respect to $K_S + 5/6 \, C$ are precisely the stable
pairs in Hacking's moduli space. All in all, we have a birational
contraction
\[
T' : \M_3 \to \HP'
\]
such that $T'(curve) = point$ if and only if $curve \subset \M_3$
is a curve in $\d_1$ obtained by varying the $j$-invariant of the
elliptic tail. Hence the forgetful map $F' : \HP' \to \M_3^{ps}$
is a bijective birational morphism between normal varieties. By
Zariski's main theorem, $F'$ is an isomorphism.
\end{proof}

\subsection{Kondo's compact moduli space} Kondo constructed a
compact moduli space of plane quartic curves by using the period
domains of K3 surfaces \cite{Kon}. Let $C$ be a smooth plane
quartic curve. Then the cyclic $\mathbb Z_4$-cover of $\P^2$
branched along $C$ is a K3 surface. The period domains of such K3
surfaces correspond to an arithmetic quotient of a bounded
symmetric domain $\mcl D$ minus two hyperplanes. He extends this
correspondence to the whole $\mcl D$ by allowing hyperelliptic
curves of genus 3 and singular pseudo-stable plane curves of genus
3 \cite{Kon}. By using the Baily-Borel compactification of period
domains, he constructs a compact moduli space which is normal and
whose boundary is one point. Details can be found in \cite{Art}
and \cite{Kon}. Let us denote Kondo's compact moduli space by
$\bar{\mcl K}$, and denote the unique point in the boundary by
$q$.

\begin{prop}\label{P:Kondo} $\bar{\mcl K}\simeq \M_3(\frac{7}{10})$.
\end{prop}

\begin{proof}
Recall from \S\ref{S:lcmodel-Chow} that $\Ch//\SL$ is isomorphic
to $ \M_3(\frac{7}{10})$. Considering Kondo's construction and the
classification of curves in $\Ch//\SL$ reveals that there is  a
forgetful morphism $F_K: \bar{\mcl K} \to \Ch//\SL$  mapping $(S,
C)$ to $C$.

We claim that there is a birational map $\HP' \to\bar{\mcl K}$
that induces an isomorphism
\begin{equation}\label{E:Kondo-isom}
\HP'\setminus Z_2\simeq \bar{\mcl K}\setminus\{ q\}.
\end{equation}
 Recall that
$\HP' = Z_0 \cup Z_1 \cup Z_2$. Let $(\P^2_{Z_0}, \mcl C)$ denote
the universal pair over $Z_0$. The universal pair has $\P^2$ as
the constant surface part and is parametrized by the curve part
that walks through all pseudostable \emph{plane} curves. By
taking the cyclic $\Z_4$-cover of $\P^2$ branched along the curve,
we obtain a universal pair $(\mcl X, \mcl C)$ of K3 surfaces
paired with pseudo-stable plane curves. This induces a morphism
from $Z_0$ to $\bar{\mcl K}$. It is an isomorphism onto its image
since it is a bijective morphism of normal varieties. A similar
construction gives an isomorphism from $Z_1$ to its image in
$\bar{\mcl K}$. Considering the description of curves in $Z_0$ and
$Z_1$, we obtain the desired isomorphism $\HP'\setminus Z_2 \simeq
\bar{\mcl K}\setminus\{ q\}$.

Retain notations from \S\ref{S:lcmodel-Chow}. The upshot of the
isomorphism (\ref{E:Kondo-isom}) is that, since $\HP'\setminus Z_2
\simeq
 \M_3^{ps}\setminus Z$ which is isomorphic under $\Psi$ to
 $\M_3^{cs}\setminus Z^{cs}$, we have
\[
\bar{\mcl K}\setminus\{ q\}\simeq (\Ch//\SL)\setminus\{
\mbox{strictly semistable point}\}.
\]
By Hartog's theorem this birational map defined away from a locus
of codimension $> 2$ is extended to an isomorphism since  both
varieties are normal.
\end{proof}

\begin{remark}  To our knowledge, the moduli space
$\Mhs$ of h-stable curves is a new modular compactification of
$M_3$. But we note that h-stable curves are precisely the curves
that appear in the semistable pairs of degree four \cite{Hac}.
Although Hacking defined the notion, the corresponding moduli
space was not constructed.
\end{remark}

\end{document}